# Nonstandard Approach to Hausdorff Measure Theory and An Analysis of Some Sets of Dimension Less Than 1

Mee Seong Im

December 2004

# University of Birmingham Research Archive

### e-theses repository

**Abstract**

We study various measure theories using the classical approach and then compute the Hausdorff dimension of some simple objects and self-similar fractals. We then develop a nonstandard approach to these measure theories and examine the Hausdorff measure in more detail. We choose to study Hausdorff measure over any other measures since it is well-defined for all sets, and widely used in many different areas in mathematics, physics, probability theory, and so forth.

Finally we generate a particular set and compute its upper and lower Hausdorff dimension. We compare our set with box-counting dimension and other well-known fractal behaviors to analyze the set in a greater detail.

# Contents

# Chapter 1

# Introduction

In this thesis we will look at some particular sets whose Hausdorff dimension is less than or equal to 1, and some of the theory of fractal sets necessary to analyze them.

We shall present some nonstandard (in the sense of Abraham Robinson [30]) techniques for defining fractal dimension. One nonstandard technique is perhaps closest to a slightly more elaborate notion of box-counting dimension which involves upper and lower limits of the number of boxes (section 6.3). We also give an alternative nonstandard method that gives exact Hausdorff dimension (section 5.2).

The sets that we will be investigating (which are defined precisely in Chapter 6) turn out to be interesting for many reasons. Our sets turn out to be versions of the Cantor set which comes up in many places and provides as a counter-example in various areas. So our particular sets may be interesting for the same reasons.

These sets seem to be new examples of fractals that are slightly more irregular than typical fractal examples. They do not resemble self-similar sets (section 6.4), and each set gives us two distinct box-counting dimensions (section 6.3). But nevertheless we can calculate their Hausdorff dimension in a satisfactory way (section 6.1).

Our sets are defined in terms of nondecreasing sequences of natural numbers and the growth rates of these sequences determine the main properties of the sets in an interesting way. We shall investigate the relationship between the Hausdorff dimension and the growth rates (section 6.1).

Finally our sets may not be self-similar, and this is interesting that we have different ways of computing the dimension. Several papers in the literature discuss the equivalence of Hausdorff dimension and box counting dimension for self-similar and other fractal sets, [10] and [11], and this work can be regarded as a contribution to this area of research.

This report is primary written for pure mathematicians with an interest in fractals, though not necessarily with any background knowledge of nonstandard analysis. The definitions and basic results of Hausdorff measure are presented from scratch. Other material that we require (such as in analysis, and measure theory in general) will be taken on trust where required.

# Chapter 2

# Basics of Nonstandard Analysis

In this chapter we will give a brief overview on the theory of nonstandard analysis. The following sources are recommended for further information regarding this chapter: [6], [20] and [29].

Before we go into any discussion we must first understand that the nonstandard approach is nothing unusual. It doesn't involve any new mathematics nor does it change the standard mathematics we're all familiar with. We have decided to use nonstandard approach because it gives us a richer meaning to standard objects and sets.

## 2.1 Nonstandard Real Numbers

**Definition 2.1.1.** *A filter over a set $I$ is an nonempty collection $\mathscr{U}$ of subsets of $I$ that satisfies:*

1. *if $A, B \in \mathscr{U}$ then $A \cap B \in \mathscr{U}$;*

2. *if $A \in \mathscr{U}$ and $B \supseteq A$, then $B \in \mathscr{U}$.*

Filters can be used to define a notion of 'large' subsets of $I$.

**Definition 2.1.2.** *A set $A \subseteq \mathbb{N}$ is considered to be large if:*

1. *a finite set including the empty set is not large;*

2. *the set of natural numbers $\mathbb{N} = \{0, 1, 2, 3, 4, \ldots\}$ is large;*

3. *if two subsets of $\mathbb{N}$ are large, then all supersets of their intersection are also large;*

4. *a set $A$ is large if and only if its complement $A^c$ is not large.*

In other words, a 'notion of largeness' is a special kind of filter called an ultrafilter satisfying the properties listed above.

**Definition 2.1.3.** *A nonprincipal ultrafilter $\mathscr{U}$ of subsets of a nonempty set $I$ contains all large sets on $I$. It is a maximal proper filter since it cannot contain any more subsets without including the empty set and becoming improper.*

We will work with nonprincipal ultrafilters which are ultrafilters with no singleton sets. And thus they have no finite sets. We note that each $I$ contains at least one maximal nonprincipal proper filter. If there is more than one nonprincipal ultrafilter, then one is isomorphic to another.

**Definition 2.1.4.** *Let $a = [a_0, a_1, \ldots]$ and $b = [b_0, b_1, \ldots]$ be sequences of real numbers. We say $a = b$ if, and only if, $a_i = b_i$ for a large set of $i \in \mathbb{N}$.*

**Remark 2.1.5.** *We will write the above definition as $(a_0, a_1, \ldots) \sim (b_0, b_1, \ldots) \Leftrightarrow \{i\colon a_i = b_i\} \in \mathscr{U}$.*

The notation of using square brackets $[a_0, a_1, \ldots]$ means that it is an equivalence class in $^*\mathbb{R}$, while the notation of using round brackets $(a_0, a_1, \ldots)$ is an element in $^*\mathbb{R}$. We will explain $^*\mathbb{R}$ later in this section.

**Lemma 2.1.6.** *For any $a = [a_0, a_1, \ldots]$, $b = [b_0, b_1, \ldots]$ and $c = [c_0, c_1, \ldots]$ the symbol $\sim$ is an equivalence relation that satisfies:*

1. $(a_0, a_1, \ldots) \sim (a_0, a_1, \ldots)$;

2. $(a_0, a_1, \ldots) \sim (b_0, b_1, \ldots) \Leftrightarrow (b_0, b_1, \ldots) \sim (a_0, a_1, \ldots)$;

3. $(a_0, a_1, \ldots) \sim (b_0, b_1, \ldots)$ and $(b_0, b_1, \ldots) \sim (c_0, c_1, \ldots)$

$$\Rightarrow (a_0, a_1, \ldots) \sim (c_0, c_1, \ldots).$$

We can easily check that the three properties are well-defined.

**Definition 2.1.7.** *Let $A$ be a set. Then the starred version of $A$, written as $^*A$, will be denoted as the nonstandard version of $A$. It is constructed by setting $^*A$ to equal $A^{\mathbb{N}}/\sim$ where $\sim$ is the equivalence relation.*

**Remark 2.1.8 (Constructing Nonstandard Real Numbers).** *We let $^*\mathbb{R} = \mathbb{R}^{\mathbb{N}}/\sim$ be the set of nonstandard reals. The algebraic operations on $^*\mathbb{R}$ are defined componentwise, and the partial order $<$ on $^*\mathbb{R}$ is defined as $a < b \Leftrightarrow \{i : a_i < b_i\} \in \mathscr{U}$. Lastly, the reals are embedded into a subset of $^*\mathbb{R}$ through the embedding $r \mapsto (r, r, r, r, \ldots)/\sim$.*

**Remark 2.1.9.** *We write $x \in {}^*\mathbb{R} \setminus \mathbb{R}$ or $x > \mathbb{R}$ to say $x$ is an infinite real. Similarly, we write $N \in {}^*\mathbb{N} \setminus \mathbb{N}$ or $N > \mathbb{N}$ to say $N$ is an infinite natural.*

## 2.2 What is an Infinitesimal?

First, we give a complete list describing the hyperreals.

**Definition 2.2.1.** *Suppose $a = [a_0, a_1, a_2, \ldots] \in {}^*\mathbb{R}$. Then for large $i \in \mathbb{N}$ we say $a$ is:*

1. *a positive infinitesimal if $0 \leq a_i < t$ for every positive real $t$;*

2. *a negative infinitesimal if $s < a_i \leq 0$ for every negative real $s$;*

3. *limited or finite if $s < a_i < t$ for some $s, t \in \mathbb{R}$;*

4. *positive unlimited if $s < a_i$ for every $s \in \mathbb{R}$;*

5. *negative unlimited if $a_i < t$ for every $t \in \mathbb{R}$;*

6. *appreciable if it is limited but is not an infinitesimal;*

7. *standard if $a_i = a_j$ for any $i, j$.*

**Example 2.2.2.** *Finite naturals $\{{}^*n\colon n \in \mathbb{N}\}$ are in the ultrafilter $\mathscr{U}$.*

**Example 2.2.3.** *If we let $a = [1/n]_{n \in \mathbb{N}}$ and $b = [1/(n^2)]_{n \in \mathbb{N}}$ be hyperreals for $n \neq 0$, then $a$ is an infinitesimal since $\{n : |\,1/n\,| < x\} \in \mathscr{U}$ for any positive real $x$. By the same reasons, $b$ is also an infinitesimal. The only difference between $a$ and $b$ is their rate of convergence. So $a > b > 0$.*

**Example 2.2.4.** *The natural $0$ is the only number that is both an infinitesimal and standard.*

## 2.3 Nonstandard Superstructure

Here we build a nonstandard superstructure from standard structures.

**Definition 2.3.1.** *Let $A$ be a set. Let $\{V_i(A)\}_{i \in \mathbb{N}}$ be the sequence defined by $V_0(A) = A$ and $V_{i+1}(A) = V_i(A) \cup \mathscr{P}(V_i(A))$. Then*

$$V(A) = \bigcup_{i \in \mathbb{N}} V_i(A)$$

*is called the superstructure over $A$. The rank of $x \in V(A)$ is the smallest $k$ so that $x \in V_k(A)$.*

We take the reader to any undergraduate's text including [12] and [13] for some introductory background on set theory and on constructing first-order logical sentences. But we will quickly remind ourselves of what consist of first-order sentences.

**Definition 2.3.2.** *First-order logic uses some of the following symbols:*

1. $\land$, $\lor$, $\neg$, $\top$, $\bot$, $=$, $\forall$, $\exists$;

2. *constant symbols such as* 0, 1;

3. *an infinite set of variables a, b, c,...;*

4. *function symbols such as* $+$, $-$, $\times$, $f$;

5. *relation symbols like* $<$, $\in$;

6. *punctuation symbols such as* $($, $)$, *and commas.*

Throughout this paper, I will use the structure that which takes the form $\mathscr{V} = \langle V, \in \rangle$. We denote $V$ as the collection of all (pure) sets, and $\in$ is the universal symbol with the meaning "the element of." The notation $\langle V, \in \rangle$ is called the first-order structure of set theory. This structure includes just about everything in the universe that one can possibly think of: all mathematics, elements, operations, etc.

We write the ultrapower of $\mathscr{V}$ as

$$^*\mathscr{V} = \prod_{\mathscr{U}} \mathscr{V} \succ \mathscr{V}.$$

This convention is convenient in terms of writing. It is easier to state, recognize and understand that our system contains all standard formalization of number theory, analysis, group theory, etc. in Zermelo and denote by the standard notation $\mathbb{R}$, $\mathbb{Q}$, $\mathbb{N}$,.... This general form also allows us to use any unexpected mathematics as they arise during our research.

Since we are going to use the universal convention $\langle V, \in \rangle$, we will not mention the structure in every chapter. In addition, we will go ahead and toss everything in $\langle V, \in \rangle$. This means we may use any standard mathematics at any time throughout this paper.

## 2.4 Łoś' Theorem

We now introduce the concept of Łoś' theorem. This idea will be used throughout our thesis.

**Theorem 2.4.1 (Łoś' Theorem).** *Let $\mathscr{V} = \langle V, \in \rangle$ be our superstructure. Let $\theta(a_0, a_1, \ldots, a_{n-1})$ be a first-order statement with $\theta({}^*a_0, {}^*a_1, \ldots, {}^*a_{n-1})$ as its $*$-transformation. Then Łoś' theorem says $\theta(a_0, a_1, \ldots, a_{n-1})$ is true in $\mathscr{V}$ if and only if $\theta({}^*a_0, {}^*a_1, \ldots, {}^*a_{n-1})$ is true in ${}^*\mathscr{V}$.*

The star that precedes an element, set, function, etc. represents that the element or the set lives in the nonstandard space, or that the operation is computed in the nonstandard world. But most of the time, we will omit the star whenever we can easily distinguish between the standard and nonstandard worlds.

When we mention $*$-transformation, we mean going from the standard world to the nonstandard world, and vice versa.

*Proof of Łoś' Theorem.* We will assume properties such as ${}^*(a+b) = {}^*a + {}^*b$, ${}^*(a \times b) = {}^*a \, {}^*\!\times {}^*b$ are preserved as we have discussed earlier. So atomic first-order statements without any connectives and quantifiers are preserved as we $*$-transfer between standard and nonstandard spaces. Now we must show that $\vee, \wedge, \neg, \rightarrow, \exists$ and $\forall$ are preserved.

We let our induction hypothesis be

$$ {}^*\mathscr{V} \models \theta({}^*a, {}^*b, {}^*c) \Leftrightarrow \{i \in \mathbb{N} : \mathscr{V} \models \theta(a_i, b_i, c_i)\} \in \mathscr{U}. \qquad (2.1) $$

This holds true since $\theta$ is a first-order statement in ${}^*\mathscr{V}$. By induction on the complexity of $\theta$, we analyze the following.

1. (the logical connective $\vee$)

    Suppose we have ${}^*\mathscr{V} \models \theta(x, y, z) \vee \psi(x, y, z)$. By the definition of $\vee$, we know ${}^*\mathscr{V} \models \theta(x, y, z)$ or ${}^*\mathscr{V} \models \psi(x, y, z)$. Since $\theta$ and $\psi$ have

complexity zero we use the induction hypothesis to obtain

$$B = \{i : \mathscr{V} \models \theta(x_i, y_i, z_i)\} \in \mathscr{U} \text{ or } C = \{i : \mathscr{V} \models \psi(x_i, y_i, z_i)\} \in \mathscr{U}.$$

Let $A = \{i : \mathscr{V} \models \theta(x_i, y_i, z_i) \vee \psi(x_i, y_i, z_i)\}$. Since $\mathscr{U}$ is an ultrafilter, $A \in \mathscr{U}$. So since $A = B \cup C$, the set $A \supseteq B$ and $A \supseteq C$.

Now to prove conversely, we have

$$A = \{i \in \mathbb{N} : \mathscr{V} \models \theta(x_i, y_i, z_i) \vee \psi(x_i, y_i, z_i)\} \in \mathscr{U}$$

for if we suppose $(B \cup C) \in \mathscr{U}$ but $B, C \notin \mathscr{U}$, then the complements $B^c, C^c \in \mathscr{U}$. If we take all possible intersections to obtain $(B \cup C) \cap B^c \cap C^c = \varnothing \in \mathscr{U}$, we have a contradiction since the empty set is in our ultrafilter. So we conclude

$$B = \{i : \mathscr{V} \models \theta(x_i, y_i, z_i)\} \in \mathscr{U} \text{ or } C = \{i : \mathscr{V} \models \psi(x_i, y_i, z_i)\} \in \mathscr{U}.$$

Then by (2.1), $^*\mathscr{V} \models \theta(x, y, z)$ or $^*\mathscr{V} \models \psi(x, y, z)$. Hence by the definition of $\vee$, $^*\mathscr{V} \models \theta(x, y, z) \vee \psi(x, y, z)$.

2. (the logical connective $\wedge$)

   The following proof holds true in both forward and in reverse direction. We have $^*\mathscr{V} \models \theta(x, y, z) \wedge \psi(x, y, z)$. By the definition of the proposition $\wedge$, $^*\mathscr{V} \models \theta(x, y, z)$ and $^*\mathscr{V} \models \psi(x, y, z)$, and by the induction hypothesis, $\{i : \mathscr{V} \models \theta(x_i, y_i, z_i)\} \in \mathscr{U}$ and $\{i : \mathscr{V} \models \psi(x_i, y_i, z_i)\} \in \mathscr{U}$. Finally, we use the intersection property of filters to obtain

   $$\{i : \mathscr{V} \models \theta(x_i, y_i, z_i) \wedge \psi(x_i, y_i, z_i)\} \in \mathscr{U}.$$

   For the reverse direction, we use the extension property of filters to obtain $\{i : \mathscr{V} \models \theta(x_i, y_i, z_i)\} \in \mathscr{U}$ and $\{i : \mathscr{V} \models \psi(x_i, y_i, z_i)\} \in \mathscr{U}$ from $\{i : \mathscr{V} \models \theta(x_i, y_i, z_i) \wedge \psi(x_i, y_i, z_i)\} \in \mathscr{U}$. Then we use all the same reasons as above to derive that $^*\mathscr{V} \models \theta(x, y, z) \wedge \psi(x, y, z)$.

3. (the connective $\neg$)

    First we begin with $^*\mathscr{V} \models \neg\, \theta(x,\, y,\, z)$. The definition of $\neg$ says $^*\mathscr{V} \not\models \theta(x,\, y,\, z)$. Using (2.1), we have $\{i : \theta(x_i,\, y_i,\, z_i)\} \notin \mathscr{U}$, and since the complement is in the ultrafilter, $\{i : \neg\, \theta(x_i,\, y_i,\, z_i)\} \in \mathscr{U}$. This proof still holds in reverse order.

4. (the "implies" symbol $\to$)

    Suppose we're given the statement $^*\mathscr{V} \models \theta(x,\, y,\, z) \to \psi(x,\, y,\, z)$. By the logical equivalence of $\to$, $^*\mathscr{V} \models \neg\, \theta(x,\, y,\, z) \vee \psi(x,\, y,\, z)$. We may now use the definition of $\neg$ and $\vee$ to obtain $^*\mathscr{V} \not\models \theta(x,\, y,\, z) \vee\, ^*\mathscr{V} \models \psi(x, y, z)$. Since the induction hypothesis tells us $\{i : \theta(x_i, y_i, z_i)\} \notin \mathscr{U}$ or $\{i : \psi(x_i,\, y_i,\, z_i)\} \in \mathscr{U}$, we can see that the complement of the first set is in the ultrafilter $\mathscr{U}$. So,

    $$B = \{i : \mathscr{V} \models \neg\theta(x_i,\, y_i,\, z_i)\} \in \mathscr{U} \text{ or } C = \{i : \mathscr{V} \models \psi(x_i,\, y_i,\, z_i)\} \in \mathscr{U}.$$

    Let $A = \{i : \mathscr{V} \models \neg\theta(x_i, y_i, z_i) \vee \psi(x_i, y_i, z_i)\}$. Since $\mathscr{U}$ is an ultrafilter, $A \in \mathscr{U}$. So since $A = B \cup C$, we may conclude $A \supseteq B$ and $A \supseteq C$.

    Conversely, we have $A = \{i : \mathscr{V} \models \neg\theta(x_i,\, y_i,\, z_i) \vee \psi(x_i,\, y_i,\, z_i)\} \in \mathscr{U}$. One of the properties of an ultrafilter says

    $$B = \{i : \mathscr{V} \models \neg\theta(x_i,\, y_i,\, z_i)\} \in \mathscr{U} \text{ or } C = \{i : \mathscr{V} \models \psi(x_i,\, y_i,\, z_i)\} \in \mathscr{U}$$

    for if we suppose $(B^c \cup C) \in \mathscr{U}$ but $B^c, C \notin \mathscr{U}$ then $B, C^c \in \mathscr{U}$. Then their intersection $(B^c \cup C) \cap B \cap C^c = (B^c \cap B \cap C^c) \cup (C \cap B \cap C^c) = \varnothing \cup \varnothing = \varnothing \in \mathscr{U}$. Because an empty set cannot be in the ultrafilter, this is a contradiction. Therefore $B^c,\, C \in \mathscr{U}$.

    Now since $B^c \in \mathscr{U}$, the set $B \notin \mathscr{U}$, and so we obtain

    $$\{i : \mathscr{V} \models \theta(x_i,\, y_i,\, z_i)\} \notin \mathscr{U} \text{ or } \{i : \psi(x_i,\, y_i,\, z_i)\} \in \mathscr{U}.$$

    We can now use (2.1) to get $^*\mathscr{V} \not\models \theta(x,\, y,\, z)$ or $^*\mathscr{V} \models \psi(x,\, y,\, z)$, and by the definition of $\neg$ and $\vee$, $^*\mathscr{V} \models \neg\, \theta(x,\, y,\, z) \vee \psi(x,\, y,\, z)$. Finally, by the propositional logic $\to$, $^*\mathscr{V} \models \theta(x,\, y,\, z) \to \psi(x,\, y,\, z)$.

5. (the existential propositional logic $\exists$)

    Suppose we have the first-order statement $^*\mathscr{V} \models \exists w \in {}^*\mathbb{R}\,\theta(x, y, z, \ldots, w)$. For some $w \in {}^*\mathbb{R}$, $^*\mathscr{V} \models \theta(x, y, z, \ldots, w)$. By our induction hypothesis, we obtain $A = \{i : \theta(x_i, y_i, z_i, \ldots, w_i)\} \in \mathscr{U}$, and so we obtain

    $$B = \{i : \exists\, w \in \mathbb{R}\, \theta(x_i, y_i, z_i, \ldots, w)\}.$$

    We conclude that $B \supseteq A$ since we cannot always determine $w_i$ from $w$. This implies $B \in \mathscr{U}$, which is what we wanted.

    Conversely let $B = \{i : \exists\, w \in \mathbb{R}\, \theta(x_i, y_i, z_i, \ldots, w)\} \in \mathscr{U}$. Then for each $i \in B$, we use the axiom of choice to choose our $w_i \in \mathbb{R}$ so that $\theta(x_i, y_i, z_i, \ldots, w_i)$ holds true. Note that if some index $i \notin B$, let $w_i = 0$ for $w = [w_0, w_1, \ldots] \in {}^*\mathbb{R}$. This process does not affect our ultrafilter since we are choosing some number that is not in $\mathscr{U}$.

    So by this definition of $w_i$, $A = \{i : \theta(x_i, y_i, z_i, \ldots w_i)\} = B \in \mathscr{U}$. By induction hypothesis, $^*\mathscr{V} \models \theta(x, y, z, \ldots w)$. Now,

    $$^*\mathscr{V} \models \theta(^*a,\, ^*b,\, ^*c) \Leftrightarrow A = \{i : \mathscr{V} \models \theta(a_i, b_i, c_i)\} \in \mathscr{U}$$

    holds true for any $^*a = [a, a, a, \ldots]$, $^*b = [b, b, b, \ldots]$ and $^*c = [c, c, c, \ldots]$. This implies $A = \mathbb{N} \cup \varnothing$. Since $\varnothing$ cannot be in the ultrafilter, $A = \mathbb{N}$.

6. (the universal propositional logic $\forall$)

    Let's suppose $^*\mathscr{V} \models \forall\, w\, \theta(x, y, z, \ldots, w)$. Then for each $w$, we have $^*\mathscr{V} \models \theta(x, y, z, \ldots, w)$. By (2.1) and for every $w = [w_0, w_1, \ldots]$, the set $A = \{i : \mathscr{V} \models \theta\,(x_i, y_i, z_i, \ldots, w_i)\}$ is in the ultrafilter. Define $B$ to be

    $$\{i : \forall\, w \in \mathbb{R}\, \mathscr{V} \models \theta(x_i, y_i, z_i, \ldots, w)\} \in \mathscr{U}.$$

    If $i \in B$, let $w_i$ be any real number. If $i \notin B$, choose $w_i \in \mathbb{R}$ so that $\mathscr{V} \models \neg\, \theta(x_i, y_i, z_i, \ldots, w_i)$. So for $w = [w_0, w_1, \ldots]$, we have

    $$A = \{i : \mathscr{V} \models \theta(x_i, y_i, \ldots, w_i)\} = B \in \mathscr{U}.$$

Now to prove conversely, we let $B$ be the set $\{i : \mathscr{V} \models \forall\, w \in \mathbb{R}\; \theta(x_i,\, y_i,\, z_i, \ldots, w)\} \in \mathscr{U}$, and let $A = \{i : \mathscr{V} \models \theta(x_i,\, y_i, \ldots, w_i)\}$. Whenever $w = w_i$ for $i \in B$, $\mathscr{V} \models \forall\, w\, \theta(x_i,\, y_i, \ldots, w)$.

Hence $A \supseteq B$, and $\mathscr{V} \models \theta(x_i,\, y_i, \ldots, w_i)$ is also true. Since the $i$'s are in $A$, $A \in \mathscr{U}$. This implies ${}^*\mathscr{V} \models \theta(x,\, y,\, z, \ldots, w)\,\forall\, w \in {}^*\mathbb{R}$. And finally, by the definition of the connective $\forall$, ${}^*\mathscr{V} \models \forall w\, \theta(x, y, z, \ldots, w)$.

This ends the proof for Łoś' theorem. □

Now let's look at an example.

**Example 2.4.2.** *Let $\theta(a)$ be the statement "$f$ is continuous at $a$." This is satisfied by Łoś' theorem since the first-order sentence*

$$\forall\, \epsilon\, [\epsilon > 0 \to \exists\, \delta\, (\delta > 0 \wedge \forall\, x\, [(\delta + x) > a\ \&\ x < (a + \delta)]$$

$$\to ((f(x) - f(a)) < \epsilon\ \&\ (f(a) - f(x)) < \epsilon))]$$

*also holds true in the nonstandard world.*

## 2.5 Internal Sets

Now we introduce the concept of internal sets.

**Definition 2.5.1.** *A subset $A$ of the nonstandard universe ${}^*V$ is called internal if there is a set $B \subseteq V$ such that*

$$A = \{[a_0,\, a_1, \ldots] \in {}^*V : \{i : a_i \in B\} \in \mathscr{U}\}.$$

In particular, $A \subseteq {}^*\mathbb{R}$ is internal if, and only if, it is ${}^*B$ for some $B \subseteq \mathbb{R}$. We note that if a set is not internal, then it is called external.

**Example 2.5.2.** *If $A = {}^*\mathbb{N} \setminus \mathbb{N}$, consider the statement: $\exists\, a\, (a = \min(A))$. This is clearly false since ${}^*B \subseteq {}^*\mathbb{N}$ has a least element. So $A$ is not internal.*

**Example 2.5.3.** *If $^*A$ is a finite subset of the hyperreals $^*\mathbb{R}$ then it is internal. For example, $\{^*0, {}^*1, {}^*2\} = {}^*\{0, 1, 2\}$. Other external sets include $\mathbb{N}$ and $\mathbb{R}$.*

Further discussion on external sets are under sections 2.6 and 2.7.

An interesting question that might arise from this is: how can we construct internal sets from standard sets? And are external sets of any interest to a mathematician? We note that we apply similar methods to construct internal elements.

Firstly to build internal sets, produce a sequence $A = \{A_i\}_{i \in \mathbb{N}}$ of standard sets so that $A = [A_0, A_1, A_2, \ldots] \in {}^*\mathscr{V}$. To generate an internal function $f$ that maps an internal set $A$ to an internal $B$, form a sequence of functions $f = \{f_i\}_{i \in \mathbb{N}}$ so that each $f_i$ is well-defined between the standard sets $A_i$ and $B_i$ for each $i$. Hence we have $f : A \to B$ if, and only if, $f_i : A_i \to B_i$ for almost all $i$. We also keep in mind that when something is internal, then there is a first-order sentence in some particular language describing the set, relation, etc.

Now to answer the second question, external sets are of no interest to a mathematician because external sets are not well-defined. When they are true in the nonstandard world, they may not be true in the standard world, and vice versa, as we have seen from previous examples. This means basic principles of mathematics are broken when moving between the standard and nonstandard worlds. So there is no great deal of interest working with external elements, sets, etc. But internal sets are again worth our attention because of its rich and consistent behavior.

We can easily see that $^*\mathbb{R}$ is much bigger than $\mathbb{R}$. Then how are the reals embedded into the hyperreals? This is where we introduce the concept of monads.

**Definition 2.5.4.** *The standard part of a finite $x \in {}^*\mathbb{R}$, written as $°x$ or*

st($x$), *is the unique real number $a$ so that it is the closest to $x$. That is,*

$$^\circ x = \text{st}(x) = \inf\{a \in \mathbb{R} : {}^*a \geq x\} = \sup\{a \in \mathbb{R} : {}^*a \leq x\}.$$

*The set of $x \in {}^*\mathbb{R}$ such that $\text{st}(x) = a$ is called the monad of $a$ and it is written as*

$$\text{st}^{-1}(a) = \mu(a) = \{x \in {}^*\mathbb{R} : \text{st}(x) = a\}.$$

We will employ the notation st($x$) to write the standard part of $x$, while both $\text{st}^{-1}(a)$ and $\mu(a)$ will be used to write the monad of $a$.

## 2.6 Returning to Infinitesimals

In this section we explore infinitesimals and their behavior with respect to both the standard and nonstandard numbers. Here, we will use Łoś' theorem throughout this section.

**Proposition 2.6.1.** *Suppose $N \in {}^*\mathbb{N} \setminus \mathbb{N}$. Then $1/N$ is an infinitesimal.*

*Proof.* Let $N \in {}^*\mathbb{N} \setminus \mathbb{N}$. Then $N > a$ for each natural number $a$. It follows that $1/N < 1/a$ for every positive $a \in \mathbb{N}$. So $1/N$ is an infinitesimal. □

**Definition 2.6.2.** *If $a \in {}^*\mathbb{R}$ is infinitesimal, we write $a \approx 0$. Similarly, if $a, b \in {}^*\mathbb{R}$, then $a \approx b$ means $a - b \approx 0$.*

**Lemma 2.6.3.** *Let $N \in {}^*\mathbb{N} \setminus \mathbb{N}$ and let $k$ be a positive integer. Then $k/N \approx 0$.*

*Proof.* We know $N > a$ for each positive $a \in \mathbb{N}$. So $1/N < 1/a$, or since $k$ is fixed, $k/N < k/a$. Since $k/a$ can be chosen as small as we like by choosing a suitable $a \in \mathbb{N}$, $k/N \approx 0$. □

In general, the above lemma holds true for a real $k$. It follows

**Remark 2.6.4.** *If $a \in {}^*\mathbb{R}$ is finite and $h \approx 0$, then $ah \approx 0$.*

**Lemma 2.6.5.** *If $a, b \in {}^*\mathbb{R}$ are finite where $b \not\approx 0$ and $h, k \approx 0$ then*

$$\frac{a+h}{b+k} \approx \frac{a}{b}.$$

*Proof.* By the previous remark, we know $\frac{h}{b+k} = \left(\frac{1}{b+k}\right) h \approx 0$. It follows

$$\frac{a+h}{b+k} = \frac{a}{b+k} + \frac{h}{b+k} \approx \frac{a}{b+k}. \tag{2.2}$$

So dividing (2.2) through by $b$, we obtain $\left(\frac{a}{b}\right) / \left(1 + \frac{k}{b}\right)$. By the continuity of the inverse map $x \mapsto x^{-1}$ (except when $x \approx 0$),

$$\left(1 + \frac{k}{b}\right)^{-1} \approx \left(1 - \frac{k}{b}\right) \approx 1.$$

So since we have approximated that

$$\frac{a/b}{(1+k/b)} \approx \frac{a/b}{1},$$

the lemma has been proved. $\square$

**Lemma 2.6.6.** *Given two finite hyperreals $a$ and $b$, $\operatorname{st}(a+b) = \operatorname{st}(a) + \operatorname{st}(b)$.*

*Proof.* Let $a_0 = \operatorname{st}(a) \in \mathbb{R}$ and $b_0 = \operatorname{st}(b) \in \mathbb{R}$. The following sentence is satisfied in the nonstandard world:

$${}^*\mathscr{V} \models \forall x \; \forall y \; \exists z \; (x = y + z).$$

Hence for $(h+k) \approx 0$, we may rewrite our $a$ and $b$ as $a = {}^*a_0 + h$ and $b = {}^*b_0 + k$. First, suppose $h > {}^*\epsilon$ where $\epsilon$ is a positive real. Then

$$\operatorname{st}(a) = \operatorname{st}({}^*a_0 + h)$$
$$= \inf\{r \in \mathbb{R} : {}^*r \geq {}^*a_0 + h\}$$
$$\geq a_0 + \epsilon.$$

This is a contradiction since $\text{st}(a) = a_0$. Similarly for some negative $\epsilon \in \mathbb{R}$, assume that $h \leq {}^*\epsilon$. Then we would find

$$\begin{aligned}\text{st}(a) &= \text{st}({}^*a_0 + h) \\ &= \inf\{r \in \mathbb{R} : {}^*r \geq {}^*a_0 + h\} \\ &\leq a_0 - |\epsilon| \\ &= a_0 + \epsilon.\end{aligned}$$

This is again a contradiction since $\text{st}(a) = a_0$. So $h$ must be an infinitesimal. Similar method is used to show $\text{st}(b) = b_0$ and that $k$ is an infinitesimal.

Now for $(h + k) \approx 0$, we observe that

$$\begin{aligned}a + b &= ({}^*a_0 + h) + ({}^*b_0 + k) \\ &= {}^*a_0 + {}^*b_0 + h + k \\ &= {}^*(a_0 + b_0) + (h + k),\end{aligned}$$

which implies that the standard part of $(a + b)$ must be $a_0 + b_0$. Hence, $\text{st}(a + b) = \text{st}(a) + \text{st}(b)$ as expected. $\square$

Finally we introduce different ways to express infinitesimals.

**Proposition 2.6.7.** *For any element $a \in {}^*\mathbb{R}$ and for some $a_0 \in \mathbb{R}$ the following are equivalent:*

1. $\text{st}(a) = a_0$;

2. $|a - {}^*a_0| \approx 0$;

3. *there is some $h \approx 0$ that satisfies $a = {}^*a_0 + h$.*

*Proof.* 1. Proving the implication from 1 to 2:

Suppose we're given $\text{st}(a) = a_0$. Then $\text{st}(a) - a_0 = 0$. If $a > {}^*a_0$ then $a - {}^*a_0 = h$ for some $0 < h \approx 0$. If $a < {}^*a_0$ then ${}^*a_0 - a = h$ for some $0 < h \approx 0$. Hence $|a - {}^*a_0| = h \approx 0$.

2. Proving the implication from 2 to 3:

   Suppose $|a - {}^*a_0| \approx 0$. Then $-k \leq a - {}^*a_0 \leq k$ for some $0 < k \approx 0$.

   Let $0 < h \approx 0$ so that $|a - {}^*a_0| = h$. If $a - {}^*a_0 = h$ then $a = {}^*a_0 + h$.

   If ${}^*a_0 - a = h$ then let $-h = h'$ so we'll have $a = {}^*a_0 + h'$.

3. Proving the implication from 3 to 1:

$$\begin{aligned}\operatorname{st}(a) &= \operatorname{st}({}^*a_0 + h) \\ &= \operatorname{st}({}^*a_0) + \operatorname{st}(h) \\ &= a_0 + 0 \\ &= a_0\end{aligned}$$

   Hence proves the proposition. $\square$

In practice, number 3 from Proposition 2.6.7 is used in almost every case since it is the most convenient mathematics when we switch back and forth from the standard and the nonstandard world.

Next, we properly show why some nonstandard sets are not internal. We begin by working with $\mu(0)$.

**Lemma 2.6.8.** *The monad of zero is not first-order definable.*

*Proof.* Let $\theta(x)$ be a first-order formula defining $\mu(0)$. Then

$$^*\mathscr{V} \models \exists\, x(x > 0 \wedge \theta(x)).$$

By Łoś' theorem, $\mathscr{V} \models \exists\, x(x > 0 \wedge \theta(x))$.

Let $\epsilon \in \mathbb{R}$ be positive so that $\theta(\epsilon)$ holds. Then $\mathscr{V} \models \theta(\epsilon)$. By Łoś' theorem again, ${}^*\mathscr{V} \models \theta({}^*\epsilon)$. Hence ${}^*\epsilon \in \mu(0)$ but this is a contradiction since ${}^*\epsilon \not\approx 0$. $\square$

**Remark 2.6.9.** *This leads us to conclude that $\mu(0)$ is not internal for if it were ${}^*A$ for some $A$, then it would be defined by*

$$\theta(x) \equiv \exists\, y(y \in {}^*A \wedge x = y).$$

**Proposition 2.6.10.** *Let $a \in \mathbb{R}$ and suppose $\mu(a)$ is the monad of $a$. Then $\mu(a)$ is not internal.*

*Proof.* Again, the proof is similar to Lemma 2.6.8. □

## 2.7 Overspill Principle

We now introduce the notion of overspill. This concept is also used repeatedly throughout this thesis. But we will mention that both Łoś' theorem and overspill principle are highly of interest whenever one approaches standard theory using nonstandard methods.

**Theorem 2.7.1 (Overspill Principle).** *Let $A \subseteq {}^*\mathbb{R}$ be a nonempty internal subset containing arbitrary large finite elements. Then $A$ contains an infinite element.*

In other words if there are arbitrarily large numbers making a first-order sentence $\theta(x)$ true, then there exists infinite natural numbers satisfying $\theta(x)$. Another perspective is if something is true for all infinitesimals, then it is true for some standard positive number. The proof for the overspill principle follows after this proposition.

**Proposition 2.7.2.** *Let $\mathbb{R}$ be a domain in $\langle \mathscr{V}, \in \rangle$. Suppose $A \subseteq {}^*\mathbb{R}$ is an internal subset and each $x \in A$ is a finite element. Then $A$ has a least upper bound.*

*Proof.* Suppose $A \subseteq {}^*\mathbb{R}$, and let $\text{st}(x) = a$. There are two cases to consider.

1. Let $A$ be bounded. So consider $\text{st}(A)$. Then the standard part of the set $A$ is also bounded since $A$ is internal. So for each $a \in \text{st}(A)$ there is $b \in \mathbb{R}$ such that $a \leq b$. So there is some smallest real $b_0$ that satisfies $a \leq b_0$ for each $a \in \text{st}(A)$. We may define $b_0 = \text{lub st}(A)$, and there is some $y \in A$ so that $y \in \text{st}^{-1}(b_0) = \text{lub}(A)$ since each $x \in A$ is finite.

2. If $A$ is unbounded then $A$ contains an arbitrary large finite $a \in {}^*\mathbb{R}$. By overspill $A$ contains an infinite element $b \in {}^*\mathbb{R}$ so that $a \leq b$ for every $a \in A$. Hence $A$ has a least upper bound.

□

*Proof of the Overspill Principle.* Let $A$ be a nonempty internal subset of ${}^*\mathbb{R}$. If $A$ is unbounded, then $A$ contains at least one infinite element and so we're done.

So we will suppose $A$ is bounded. Let $a$ be the least upper bound of $A$. Since $A$ contains arbitrary large finite elements, $a$ must be infinite. If there is no $x \in A$ so that $a - \epsilon \leq x \leq a$ for some positive $\epsilon$, then $a - \epsilon$ is the least upper bound of $A$, which is a contradiction.

So there is some $x \in A$ so that $a - \epsilon \leq x \leq a$, and hence proves the theorem. □

**Remark 2.7.3.** *Overspill may also be used to show other sets, including the reals, are not internal to ${}^*\mathscr{V}$ for if they are internal, then we would be assuming that $\mathbb{R} = {}^*\mathbb{R}$.*

We refer to [24] to see more examples of external sets, and to see some applications of nonstandard analysis, we take the reader to [19].

Finally, we introduce the concept of saturation. This concept will not be used too often in this thesis since the earliest we will apply the theory of saturation will be in the proof for Theorem 5.2.2.

**Theorem 2.7.4 (Saturation).** *Let $\{A^i\}_{i \in \mathbb{N}}$ is a sequence of internal sets that satisfies $\bigcap_{i=0}^{n} A^i \neq \varnothing$ for each $n \in \mathbb{N}$ where each $A^i = [A_0^i, A_1^i, A_2^i, \ldots]$. Then $\bigcap_{i \in \mathbb{N}} A^i \neq \varnothing$.*

In "Nonstandard Analysis and its Applications," N. Cutland calls this saturation as $\aleph_1$-saturation. Here, we will leave it as it is. We take the reader to [6] for the proof of Theorem 2.7.4.

# Chapter 3

# Measure Theory Review

Here we will mention some of history's mathematical analysts and briefly share their personal lives. We will also discuss their contribution in the area of analysis and measure theory.

## 3.1 History

We begin by mentioning a German mathematician Georg Friedrich Bernhard Riemann, born on September 17, 1826. He made an important contribution in analysis and in differential geometry. Riemann was the second child out of his six siblings. His father was a Lutheran minister who had a great desire for Riemann to study theology. But Bernhard showed particular interest in mathematics since his youth. He studied mathematical texts from his secondary school where the director of the school gave Bernhard an access to his library. It is believed that Bernhard borrowed a 900-page book on number theory by Adrien-Marie Legendre that which he read in less than a week.

When Riemann was 19 years of age he enrolled at the University of Göttingen, Germany to study theology. But with the permission of his father, he soon transfered to the department of philosophy in order to further his knowledge in mathematics. Riemann's gift in mathematics wasn't discovered until he transfered to Berlin University in 1847 to work under Ferdinand Gotthold Max Eisenstein and Johann Peter Gustav Lejeune Dirichlet. Also through Johann Benedict Listing and Wilhelm Eduard Weber, Riemann heavily studied theoretical physics and topology [26].

In 1849, he returned to the University of Göttingen to work on his Ph. D. under Johann Carl Friedrich Gauss. His thesis consisted of what is known today as Riemann surfaces. It composed of geometric properties of conformal mappings which are surface maps, such as $\mathbb{R}^2 \to \mathbb{R}^2$ or $\mathbb{C} \to \mathbb{C}$, which preserve angles, surface connectivity and analytic functions that are differentiable at every point in the given domain.

After obtaining his doctorate degree Riemann contributed toward many branches of mathematics such as geometry, abelian functions, topology and number theory. In this paper we will concentrate on his work in analysis. He came up with the idea known as Riemann integral. The basic idea is the following: let $[a,\ b] \subseteq \mathbb{R}$ and suppose $f$ maps $[a,\ b] \to \mathbb{R}_+$. Then for each $x$ in the domain, $|f(x)| \leq M$ where $M$ is a positive real. Then Riemann integral is the area under the function $f$ as $x$ varies from $a$ to $b$.

In other words, we partition the interval $[a,\ b]$ into finitely-many subintervals $[x_{i-1},\ x_i)$ where $a = x_0 < x_1 < x_2 < \ldots < x_n = b$, and we let $s$ to be the function

$$s(\ [x_{i-1},\ x_i)\ ) = f(x_{i-1})$$

for each $i = 1, \ldots n$. Then

$$\sum_{i=1}^{n} s\left([x_{i-1},\ x_i)\right) |x_i - x_{i-1}|$$

is called the step function. If $f$ is monotonically nondecreasing and positive then $s$ is bounded by

$$0 \leq \sum_{i=1}^{n} s\left([x_{i-1},\ x_i)\right) |x_i - x_{i-1}| \leq \int_a^b f(x)\, dx.$$

But if not, the function $s$ over $[a, b]$ is bounded by some positive constant $M$ times the width of the interval $[a, b]$ [5].

A modification on the Riemann integral uses the notion of upper and lower sums. Let $f\colon [a, b] \to \mathbb{R}$ be a function. Partition the interval so that

$$[a, b] \supseteq \bigcup_{i=1}^{n}(x_{i-1}, x_i)$$

where $a = x_0 < x_1 < x_2 < \ldots < x_n = b$. We will define $\overline{f}(x_{i-1})$ to be $\sup\{f(x)\colon x \in (x_{i-1}, x_i)\}$ and $\underline{f}(x_{i-1}) = \inf\{f(x)\colon x \in (x_{i-1}, x_i)\}$. Then the Riemann upper sum of $f$ is said to be

$$\sum_{i=1}^{n} \overline{f}(x_{i-1}) \left|x_i - x_{i-1}\right|,$$

and the Riemann lower sum of $f$ is

$$\sum_{i=1}^{n} \underline{f}(x_{i-1}) \left|x_i - x_{i-1}\right|.$$

It follows that

$$\sum_{i=1}^{n} \underline{f}(x_{i-1}) \left|x_i - x_{i-1}\right| \leq \sum_{i=1}^{n} \overline{f}(x_{i-1}) \left|x_i - x_{i-1}\right|.$$

Now if we take $n$ to be sufficiently big so that

$$\lim_{n \to \infty} \sum_{i=1}^{n} \underline{f}(x_{i-1}) \left|x_i - x_{i-1}\right| = \lim_{n \to \infty} \sum_{i=1}^{n} \overline{f}(x_{i-1}) \left|x_i - x_{i-1}\right|,$$

then the integral $\int f(x)dx$ over $a$ to $b$ exists and is called the Riemann sum [5]. We also say $f$ is Riemann integrable whenever the Riemann upper sum of $f$ equals the Riemann lower sum of $f$. In general, the integral of a function has the following properties: for any function $f$,

1. absolute: $\left|\int_A f\right| \leq \int_A |f|$ for any $A \subseteq \mathbb{R}$;

2. additivity: $\int_{A \cup B} f = \int_A f + \int_B f$ where $A \cap B = \varnothing$ for any nonempty sets $A, B \subseteq \mathbb{R}$;

3. linearity: $\int_A (af + bg) = a \int_A f + b \int_A g$ for any $a, b \in \mathbb{R}$;

4. monotone: $f(x) \leq g(x)$ for each $x \in A$ implies $\int_A f \leq \int_A g$.

In 1863 about a year after his marriage to Elise Koch, Riemann caught a cold which led to tuberculosis. In fact, the Riemann family had a history of health problems, which led for him to lose his mother and most of his siblings at a tender age. He tried to win against this illness by temporary relocating to Italy, which is warmer than Germany. Yet, his health continued to fail rapidly, and on July 20, 1866 he died in Selasca, Italy. His lack of rigorous proofs led many to doubt his work. But through Hermann Amandus Schwarz, Felix Christian Klein, David Hilbert and others, Riemann's approach was mended, producing important results for both physicists and mathematicians.

An improvement of Riemann integral was made when mathematicians came across functions that could not be integrated simply by using Riemann's methods. An example of such a function is the characteristic function $\chi_A(x)$. We let $A$ be the set of rationals in [0, 1]. The characteristic function gives 1 if $x \in \mathbb{Q}$, and it equals 0 whenever $x \notin \mathbb{Q}$. This function is not Riemann integrable. So how shall we go about integrating $\chi_A(x)$ over the unit interval?

Henri Léon Lebesgue was born on June 28, 1875 in Beauvais, France. He first studied at the Beauvais College and then traveled to Paris for higher education. In 1894 he became a student at the École Normale Supérieure in Paris, and then he received his teaching diploma in mathematics three years later. His next two years was spent reading René-Louis Baire's papers at the École Normale Supérieure library. In 1901, Lebesgue first proposed measure theory and by the following year, he introduced Lebesgue integral which extended the Riemann integral to a wide variety of discontinuous functions. His idea is discussed in section 3.3.

Lebesgue received many awards, including the Prix Houllevigue (1912) and the Prix Saintour (1917). He was appointed to various university chairs in France and wrote approximately 90 papers and books which include areas in potential theory, the calculus of variations, set theory and dimension theory. The last 20 years of his life was dedicated towards teachers and education as he worked on geometry and history of mathematics. Lebesgue died on July 26, 1941 in Paris, France [27].

Then came Felix Hausdorff who was born on November 8, 1868 in Breslau, Germany. He graduated from Leipzig at the age of 23 and later became a professor there for about 20 years. Although he was offered a position at the University of Göttingen, he continued working at Bonn from 1910 to 1935.

Most of Hausdorff's contributions were in the areas of topology and set theory. He studied partial order relation, and he introduced ordinals to attempt to prove Georg Cantor's continuum hypothesis, which says there is no set with cardinality that is strictly between the naturals $\mathbb{N}$ and the reals $\mathbb{R}$. He also asked if $2^{\aleph_n} = \aleph_{n+1}$. That is, for any set $A$ there is no set $B$ so that $\text{card}(A) \lneqq \text{card}(B) \lneqq \text{card}(\mathscr{P}(A))$. He also developed what is known today as Hausdorff dimension. The Hausdorff dimension of a set turns out to be greater than or equal to the topological dimension of a set and less than or equal to the dimension of the metric space imbedding the set. In addition, he concluded that it does not need to be an integer.

In 1935 Hausdorff was forced to retire due to the Nazi Nuremberg Laws. But he continued studying mathematics, and his research could only be published outside Germany. In 1941, he had to go to an internment camp because of his Jewish background. He was able to avoid it until the following year when he could no longer avoid the camp. On January 26, 1942 Hausdorff, his wife and his wife's sister committed suicide in Bonn, Germany [25], [38].

Although his life was ended by inevitable circumstances, he is one of the distinguished mathematicians of all time.

## 3.2 Introducing the Concept of Measure

First we will introduce what we mean by measure.

**Definition 3.2.1.** *Let $\mathscr{A}$ be a set of all internal subsets of $X$. A measure is a nonnegative function $\mu : \mathscr{A} \to \mathbb{R}_+ \cup \{0\}$ satisfying:*

1. *$\mu(\varnothing) = 0$;*

2. *$\mu\left(\bigcup_i A_i\right) = \sum_i \mu(A_i)$ for $A_i \cap A_j = \varnothing$ for any $i \neq j$.*

A measure on a set $A \subseteq X$ is similar to the concept of "the size of $A$" which is usually denoted as $\mu$ or $\mu(A)$. It determines lengths, areas, volumes, probabilities, etc. to a given set.

**Definition 3.2.2.** *A $\sigma$-algebra, or $\sigma$-field, over a set $X$ is the power set of $X$ denoted as $\mathscr{P}(X)$ which is closed under countable set operations. We have $\sigma$-algebra satisfying the following:*

1. *the empty set $\varnothing$ is in $\sigma$-algebra;*

2. *if $A \subseteq X$ is in $\sigma$-algebra then $X \setminus A$ is in $\sigma$-algebra;*

3. *if $A_0, A_1, A_2, \ldots$ is a sequence in $\sigma$-algebra then $\bigcup_{i \in \mathbb{N}} A_i$ is in $\sigma$-algebra.*

The sets in $\sigma$-algebra are called measurable sets. Moreover, a set $X$, together with the collection of measurable sets and measure $\mu$, is called a measure space.

**Lemma 3.2.3.** *The measure $\mu$ has the following properties:*

1. *$\mu(A) \geq 0$ for any $A \subseteq X$;*

2. *if $\{A_i\}_{i=0}^m$ is a collection of disjoint sets, then*

$$\mu\left(\bigcup_{i=0}^m A_i\right) = \sum_{i=0}^m \mu(A_i).$$

In addition if $A_0, A_1, A_2, \ldots A_m, \ldots$ is a sequence of disjoint sets then

$$\mu\left(\bigcup_{i \in \mathbb{N}} A_i\right) = \sum_{i \in \mathbb{N}} \mu(A_i);$$

3. if $A \cong B$ ($A$ is congruent to $B$) then $\mu(A) = \mu(B)$.

**Definition 3.2.4.** *If $A \subseteq X$ is a measurable set, then $A$ is called a null set if $\mu(A) = 0$. The measure $\mu$ is called complete if for any subset $B$ of a null set $A$, $B$ is measurable. Moreover $B$ is also a null set.*

Define $V(A)$ to be the product of the lengths of the involved intervals of $A$. We say a set $A \subseteq \mathbb{R}^n$ is a null set if for every $\epsilon > 0$, the set $A$ can be written as $\bigcup_{i \in \mathbb{N}} U_i$ where $V(\bigcup_{i \in \mathbb{N}} U_i) \leq \epsilon$.

**Proposition 3.2.5.** *Let $A_0, A_1, A_2, A_3, \ldots$ be a countable sequence of measurable sets. The properties for the measure $\mu$ on $\sigma$-algebra include:*

1. $\mu(A_0) \leq \mu(A_1)$ if $A_0 \subseteq A_1$;

2. $A = \bigcup_{i \in \mathbb{N}} A_i$ is measurable and $\mu(A) = \lim_{i \to \infty} \mu(A_i)$ if $A_i \subseteq A_{i+1}$ for every $i$;

3. $A = \bigcap_i A_i$ is measurable if $A_i \supseteq A_{i+1}$ for each $i$; and if there is some $A_i$ so that $\mu(A_i) < \infty$, then $\mu(A) = \lim_{i \to \infty} \mu(A_i)$.

**Definition 3.2.6.** *Let $\mu : \mathscr{P}(\mathbb{Z}) \to [0, 1]$ satisfy:*

1. $\mu(\mathbb{Z}) = 1$;

2. $\mu(A \cup B) = \mu(A) + \mu(B)$ for $A, B \subseteq \mathbb{Z}$ and $A \cap B = \emptyset$;

3. $\mu(A_n) = \mu(A)$ where $A_n = \{a + n : a \in A\}$ for any $n \in \mathbb{Z}$.

*Then the map $\mu$ is called an invariant mean.*

Moreover if $C = A \cup B$ with $A \cap B = \emptyset$, then we have

$$\mu(A) \leq \mu(A) + \mu(B) = \mu(A \cup B) = \mu(C).$$

So $\mu$ is monotonically increasing map.

## 3.3 Lebesgue Measure

Now, since we have seen some definitions and properties of measure, we begin by briefly introducing Lebesgue measure.

**Definition 3.3.1 (Lebesgue Measure).** *Let $A$ be the union of the intervals $(a_i, b_i) \subseteq \mathbb{R}$ where $(a_i, b_i) \cap (a_j, b_j) = \varnothing$ for $i \neq j$. Then we define its Lebesgue measure as*

$$\lambda(A) = \sum_i (b_i - a_i).$$

*In general, if*

$$A = \bigcup_j \left[ \prod_{i=1}^{k_j} (a_i, b_i) \right]_j \subseteq \mathbb{R}^n$$

*and the sets $\prod_i (a_i, b_i)$ are disjoint where each $1 \leq k_j \leq n$ then*

$$\lambda(A) = \sum_j \left[ \prod_{i=1}^{k_j} (b_i - a_i) \right]_j.$$

Hence the Lebesgue measure of $A$ is the product of intervals for which $A$ is defined. Lebesgue measurable sets form a $\sigma$-algebra and $\lambda$ is the unique translation invariant measure on that $\sigma$-algebra. Now, we will introduce the definition of Lebesgue outer measure.

**Definition 3.3.2.** *Let $A$ be a subset of $\mathbb{R}^n$ and suppose $U_j$ is a clopen (both closed and open) subset of $\mathbb{R}^n$ that satisfies*

$$\bigcup_j U_j = \bigcup_j \left[ \prod_{i=1}^{k_j} [x_i, y_i) \right]_j \supseteq A$$

*where each $[x_i, y_i) \subseteq \mathbb{R}$ and $x_i < y_i$. The volume of each $U_j$ is computed as*

$$V(U_j) = \left[ \prod_{i=1}^{k_j} (y_i - x_i) \right]_j.$$

*Then Lebesgue $n$-dimensional outer measure is defined to be*

$$\overline{\lambda}(A) = \inf \sum_j V(U_j)$$

*where the collection of $U_j$ covers the set $A$.*

**Definition 3.3.3.** *Let $A \subseteq X$, and let $a > 0$. Lebesgue inner measure of $A$ is defined to be*

$$\underline{\lambda}(A) = a - \overline{\lambda}(X \setminus A)$$

*where $\overline{\lambda}(X \setminus A)$ is the Lebesgue outer measure of $X \setminus A$.*

Note that if $X = [0, 1]$, then $a = 1$.

**Proposition 3.3.4.** *A subset $A \subseteq X$ is measurable if for any $B \subseteq X$,*

$$\overline{\lambda}(A) = \underline{\lambda}(A) \Leftrightarrow \overline{\lambda}(B) = \overline{\lambda}(B \cap A) + \overline{\lambda}(B \setminus A).$$

The proof is found in B. Craven's text [5]. Note that if $\overline{\lambda}(A) = \underline{\lambda}(A)$, then the number obtained is the Lebesgue measure of $A$, which is normally written as $\lambda(A)$ [8].

**Proposition 3.3.5.** *Lebesgue measure has the following properties:*

1. *a unit interval $[0, 1]$ has Lebesgue measure 1;*

2. *$\lambda([0, 1] \times [0, 1] \times \ldots \times [0, 1]) = 1$, and if*

$$A = \prod_{i=1}^{k} I_i \subseteq \mathbb{R}^n$$

   *where each $I_i \subseteq \mathbb{R}$ is an interval then*

$$\lambda(A) = |I_1| \times |I_2| \times \ldots \times |I_k|$$

   *where $|I_i|$ is the length of $I_i$;*

3. *Lebesgue measure holds all $\sigma$-algebra and measure $\mu$ properties;*

4. *if $A$ is a Lebesgue measurable null set then for any $B \subseteq A$, $B$ is also a Lebesgue measurable null set;*

5. *if $A \subseteq \mathbb{R}^n$ is Lebesgue measurable and $x \in \mathbb{R}^n$ then*

$$A + x = \{a + x : a \in A\}$$

   *is also Lebesgue measurable, and $\lambda(A) = \lambda(A + x)$. This is called the translation of $A$ by $x$.*

We do not need to discuss Lebesgue measure in too much detail since it is defined only for some sets. It cannot be extended to every subsets of $\mathbb{R}^n$ because Lebesgue measure does not remain unchanged under translation and rotation. See [34] for an example.

We refer to [4], [5] and [8] for more studies on measure theory. We will now move-on to Hausdorff measure.

## 3.4 Hausdorff Measure

Let us introduce few classical definitions from [8] and [31].

**Remark 3.4.1.** *The function d: $\mathbb{R}^n \times \mathbb{R}^n \to \mathbb{R}$ takes any $x$, $y \in \mathbb{R}^n$ and gives a real value*

$$d(x, y) = ||x - y||.$$

*For any $x$, $y$, $z \in \mathbb{R}$, a metric is a function $d : \mathbb{R}^2 \to \mathbb{R}_+ \cup \{0\}$ which satisfies:*

1. *$d(x, y) \geq 0$, or equals 0 if $x = y$;*

2. *$d(x, y) = d(y, x)$;*

3. *$d(x, y) \leq d(x, z) + d(z, y)$.*

The notation $d(x, y)$ gives us the flexibility to choose whichever metric we happened to be using. We will write our metric as $d(x, y)$ throughout this paper but we will assume, more specifically, the standard Euclidean metric

$$d(x, y) = \sqrt{\sum_{i=1}^{n}(x_i - y_i)^2}$$

throughout our analysis.

**Definition 3.4.2.** *We define the diameter of $U \subseteq \mathbb{R}^n$ as*

$$\mathrm{diam}(U) = \sup\{d(x, y) : x, y \in U\}$$

*where d is a metric.*

**Definition 3.4.3.** *Let $A \subseteq \mathbb{R}^n$ be in a metric space. For any $U_i \subseteq \mathbb{R}^n$, we say the collection $\{U_i : i \in \mathbb{N}\}$ is a $\delta$-cover of A if:*

1. *$\bigcup_{i \in \mathbb{N}} U_i \supseteq A$;*

2. *$0 \leq \operatorname{diam}(U_i) \leq \delta$ for each $i \in \mathbb{N}$.*

We note that we haven't described how the covering $\{U_i\}$ of $A$ looks like. The $\delta$-covers may be rectangular, circular, triangular, spherical or some odd shape. The only condition we care about is the length of the diameter of each covering. We will let the reader decide on the shape of the covering, although choosing the shape actually depends on $A$. It should be chosen so that the difference $\bigcup U_i \setminus A$ and $A \setminus \bigcup U_i$ is a negligible set. As we have mentioned under Definition 3.2.4, a negligible set is a set that is small enough so that it can be ignored for our purposes. Negligible sets form an ideal: the empty set is negligible, the union of two negligible sets is negligible, and any subset of a negligible set is negligible. Moreover, we will assume that it is sigma-ideal: countable unions of negligible sets are also negligible.

Now let us introduce the definition of Hausdorff measure.

**Definition 3.4.4.** *Suppose $A$ is a subset of $\mathbb{R}^n$. Let $s$ and $\delta$ be positive. Then the Hausdorff $s$, $\delta$-measure, or simply Hausdorff measure, is defined to be*

$$H_\delta^s(A) = \inf \left\{ \sum_{i=0}^{\infty} \operatorname{diam}(U_i)^s \right\}$$

*where the infimum is over all $\delta$-covers $\{U_i\}_{i \in \mathbb{N}}$ of $A$.*

We will also assume that the $U_i$'s which we will encounter repeatedly are measurable. Similar to Proposition 3.3.4, a classical method to determine the measurability of a set $A \subseteq \mathbb{R}^n$ is whenever, given an outer measure $\mu$ on $\mathbb{R}^n$ and for any $B \subseteq \mathbb{R}^n$,

$$\mu(B) = \mu(B \cap A) + \mu(B \cap A^c)$$

holds, where $A^c = \mathbb{R}^n \setminus A$. We will now introduce the Hausdorff outer measure.

**Definition 3.4.5.** *Given a set $A \subseteq \mathbb{R}^n$ and $s > 0$ the Hausdorff $s$-dimensional outer measure of $A$ is*

$$H^s(A) = \lim_{\delta \to 0} H^s_\delta(A) = \sup_{\delta > 0} H^s_\delta(A).$$

Now we discuss $H^s(A)$ in greater detail. We begin by explaining why the Hausdorff $s$-dimensional outer measure of $A$ is the supremum of all such $H^s_\delta(A)$'s.

**Remark 3.4.6.** *The Hausdorff measure $H^s_\delta$ increases as $\delta$ decreases.*

*Proof.* Let $0 < \epsilon < \delta$. Suppose the collection $\{U_i\}$ is an $\epsilon$-cover of $A$. Since $0 < \text{diam}(U_i) \leq \epsilon < \delta$, it is also a $\delta$-cover of $A$ by definition. Then we have

$$H^s_\delta(A),\ H^s_\epsilon(A) \leq \sum_i \text{diam}(U_i)^s$$

since both $H^s_\delta(A)$ and $H^s_\epsilon(A)$ take infimum over all possible $\delta$, $\epsilon$-cover of $A$ respectively. Now taking the infimum over $\epsilon$-covers, we obtain

$$H^s_\delta(A) \leq \inf_{\epsilon > 0} \sum_i \text{diam}(U_i)^s = H^s_\epsilon(A).$$

Hence we want $H^s(A)$ to be the supremum of $H^s_\delta(A)$'s as $\delta \to \epsilon$. □

Next, we observe that as $s$ ranges over 0 to infinity, Hausdorff $s$-dimensional outer measure $H^s(A)$ is a nonincreasing function.

**Remark 3.4.7.** *Let $A \subseteq \mathbb{R}^n$. If $0 \leq s < t$ then $H^s(A) \geq H^t(A)$.*

So if some number $s$ approaches the correct Hausdorff dimension $\dim_H(A)$ from above, the Hausdorff $s$-dimensional outer measure will increase. If it's approaching the correct $\dim_H(A)$ from below, then $H^s(A)$ will decrease. Since the limit of $H^s(A)$ always exist [17], the critical point $s$ always exist.

*Proof.* Let $0 < \delta < 1$ and the collection $\{U_i\}$ is a $\delta$-cover of $A$. Then we have $\bigcup_i U_i \supseteq A$ and each $\text{diam}(U_i) < 1$. We observe that for each $i$,

$$\text{diam}(U_i)^t \leq \text{diam}(U_i)^s$$

and $\delta^s \leq \delta^t$. Since the following inequality holds true

$$\delta^s \times H^t_\delta(A) = \delta^s \left( \inf \sum_i \text{diam}(U_i)^t \right)$$
$$\leq \delta^t \left( \inf \sum_i \text{diam}(U_i)^s \right)$$
$$= \delta^t \times H^s_\delta(A),$$

we conclude

$$\delta^{s-t} H^t_\delta(A) \leq H^s_\delta(A),$$

where $\delta^{s-t}$ is a positive real. Now if $H^t_\delta(A)$ is finite and strictly positive, $H^s_\delta(A)$ tends to infinity. So we observe that if $0 < H^t(A) < \infty$, then the critical point is at $t$. And since $s < t$, $H^s(A) \to \infty$. Hence $H^t(A) \leq H^s(A)$, and Hausdorff $s$-dimensional outer measure is a nonincreasing function. □

**Definition 3.4.8.** *Let $A$ be a subset of $\mathbb{R}^n$. The Hausdorff dimension of $A$ is defined as*

$$\dim_H(A) = \inf\{s : H^s(A) = 0\} = \sup\{s : H^s(A) = \infty\}.$$

If $s$ is greater than the actual Hausdorff dimension of $A$, then we see that $H^s(A) = 0$. So we define the Hausdorff dimension of $A$ to be $\dim_H(A)$ = $\inf\{s\colon H^s(A) = 0\}$. And if $s$ is less than the actual $\dim_H(A)$, then $H^s(A) = \infty$, which gives us $\dim_H(A) = \sup \{s\colon H^s(A) = \infty\}$.

### 3.4.1 Hausdorff Density

Before we discuss Hausdorff density, let's first recall the definition of the classical Lebesgue density. It says for any nonempty $A \subseteq [0, 1]$ and its

$\lambda$-measure, the density of $A$ at $x$ is

$$d_\lambda(A,\ x) = \lim_{h \to 0^+} \frac{\lambda(A \cap [x-h,\ x+h])}{\lambda([x-h,\ x+h])}.$$

So a set is considered regular if, and only if, $d_\lambda(A,\ x) = 1$ for almost all $x \in A$. If not, then $A$ is said to be irregular.

K. J. Falconer's "The Geometry of Fractal Sets" [8] uses similar definition to obtain the density of a set at a point $x$. First we let $B_r(x)$ be the (open/closed) ball of radius $r$ centered at $x$ so that $\mathrm{diam}(B_r(x)) = 2r$.

**Definition 3.4.9.** *Let $A \subseteq \mathbb{R}^n$ be a nonempty set with Hausdorff dimension $s$. Then the lower and upper densities of $A$ at a point $x \in \mathbb{R}^n$ are defined as*

$$\underline{D}_s(A,\ x) = \liminf_{r \to 0} \frac{H^s(A \cap B_r(x))}{\mathrm{diam}(B_r(x))^s}$$

*and*

$$\overline{D}_s(A,\ x) = \limsup_{r \to 0} \frac{H^s(A \cap B_r(x))}{\mathrm{diam}(B_r(x))^s}$$

*respectively. If $\underline{D}_s(A, x) = \overline{D}_s(A, x)$ then we say the density of $A$ at $x$ exists and we write $D_s(A,\ x)$.*

The Hausdorff $s$-dimensional outer measure approximates the expected measure of $A$ while $D_s(A,\ x)$ gives a local measure of our set $A$. Similar to the classical Lebesgue density, $D_s(A,\ x) = 1 \Leftrightarrow x \in A$ is regular, and $D_s(A,\ x) \neq 1 \Leftrightarrow x \in A$ is irregular. The rest of K. J. Falconer's book [8] characterizes regular sets and obtain bounds for the densities of various sets. All of these eventually lead us to prove that a set cannot be regular unless its Hausdorff dimension $s$ is an integer.

## 3.5 Box-counting Measure

Now let's consider the concept behind the box-counting dimension. Suppose we want to cover a line of length 1 with small line segments of length $\epsilon > 0$. Then how many will we need? It's obvious that we need $1/\epsilon$ of them to cover the line. Now suppose we have a unit square and we want to cover it with smaller squares of side-length $\epsilon$. Then we would need $1/\left(\epsilon^2\right)$ of them. And as to cover a unit cube with smaller cubes of side-length $\epsilon$, we would need $1/\left(\epsilon^3\right)$ of them. So the box-counting dimension goes as follows [14], [33]:

**Definition 3.5.1.** *Let $\epsilon > 0$, and let $A \subseteq \mathbb{R}^n$. Let $N_\epsilon(A)$ denote the number of n-dimensional cubes with side-length $\epsilon$. If there is some positive number $d_{bc}$ so that the ratio of $N_\epsilon(A)$ over $1/\left(\epsilon^{d_{bc}}\right)$ is a positive constant scale-factor as $\epsilon \to 0$, then we say $d_{bc}$ is the box-counting dimension of A.*

In other words for some positive $k \in \mathbb{R}$ if

$$\lim_{\epsilon \to 0} \frac{N_\epsilon(A)}{1/\left(\epsilon^{d_{bc}}\right)} = k, \tag{3.1}$$

then $d_{bc} = \dim_{bc}(A)$ is the box-counting dimension of $A$. Let's consider few examples.

**Example 3.5.2.** *Let A be a unit square. Given smaller squares of side-length $\epsilon > 0$ we need exactly $1/\left(\epsilon^2\right)$ of them to cover A as we have said above. So we will let $N_\epsilon(A) = 1/\left(\epsilon^2\right)$. Suppose we let $d_{bc} = 1$. Then according to the limit equation above, we have*

$$\lim_{\epsilon \to 0} \frac{1/\left(\epsilon^2\right)}{1/\epsilon} = \lim_{\epsilon \to 0} \frac{1}{\epsilon} = \infty.$$

*So $d_{bc}$ is not equal to 1. Suppose we let $d_{bc} = 3$. Then similar as above,*

$$\lim_{\epsilon \to 0} \frac{1/\left(\epsilon^2\right)}{1/\left(\epsilon^3\right)} = \lim_{\epsilon \to 0} \epsilon = 0$$

*Hence in order for the above limit to converge to some positive real, we can easily see that the box-counting dimension must be 2.*

The relation between $N_\epsilon(A)$ and $1/\left(\epsilon^{d_{bc}}\right)$ is clear when $d_{bc}$ is equal to the topological dimension. But for more complicated sets, $d_{bc}$ may not be so clear.

Let's go back to equation (3.1) and take the log of both sides to obtain

$$\lim_{\epsilon \to 0} (\log N_\epsilon(A) + d_{bc} \log \epsilon) = \log k.$$

We solve for $d_{bc}$ to get

$$d_{bc} = \lim_{\epsilon \to 0} \frac{\log k - \log N_\epsilon(A)}{\log \epsilon}.$$

Since the term $\log k$ is a fixed real number and $(\log \epsilon) \to \infty$ as $\epsilon \to 0$, the first term drops out leaving us with

$$d_{bc} = \lim_{\epsilon \to 0} \frac{-\log N_\epsilon(A)}{\log \epsilon}.$$

The signs correspond with the sign of $d_{bc}$ since for each small $\epsilon$, we observe that $\log \epsilon < 0$.

We refer to [9], [22] and [37] for the next definition.

**Definition 3.5.3.** *Let $\epsilon > 0$, and let $A$ be a subset of $\mathbb{R}^n$. Suppose $N_\epsilon(A)$ is the number of n-dimensional cubes of volume $\epsilon^n$ that covers $A$. Then the lower box-counting dimension is defined as*

$$\underline{d}_{bc}(A) = \liminf_{\epsilon \to 0} \frac{-\log N_\epsilon(A)}{\log \epsilon}$$

*and the upper box-counting dimension is*

$$\overline{d}_{bc}(A) = \limsup_{\epsilon \to 0} \frac{-\log N_\epsilon(A)}{\log \epsilon}.$$

*If $\underline{d}_{bc}(A) = \overline{d}_{bc}(A)$ then the limit is said to exist and equals the box-counting dimension $d_{bc}(A)$.*

Since box-counting measure uses sets of equal diameter as opposed to the Hausdorff measure (which takes the infimum over all sets), the box-counting theorem is much easier to calculate. But when it comes to calculating sets of irregular or discontinuous shapes and patterns, then the box-counting theorem fails in the accuracy of finding the correct dimension that stays consistent despite the variation of the sets [22].

Today the box-counting dimension is used to classify various textures, patterns and surfaces in an image. And although the true dimension of the set may be incorrect up to 40% when using the box-counting definition [23], $d_{bc}$ is a dimension in its own right which may be used to determine the complexity of a set. D. Neary discussed that the box-counting methods perform badly for sets which do not have universal self-similarity [23]. But for more regular sets (like a walnut with less directional surface), the box-counting dimension correspond with the rest of its computation on roughness. Hence $d_{bc}$ should be considered for more predictable sets, but for irregular sets, the box-counting dimension results are less dependable.

Finally for any set $A$, we note that the topological dimension $\dim_T(A)$ is less than the Hausdorff dimension $\dim_H(A)$ which is again less than the box-counting dimension [32], i.e.,

$$\dim_T(A) \leq \dim_H(A) \leq \underline{\dim}_{bc}(A) \leq \overline{\dim}_{bc}(A).$$

This tells us that $\dim_{bc}(A)$ normally overestimates, and that the box-counting dimension may give us distinct lower and upper dimensions, which will depend on $A$.

Lastly we introduce a property regarding box-counting dimension.

**Proposition 3.5.4.** *Let $A$ be a nonempty subset of $\mathbb{R}^n$. Then*

$$\dim_{bc}(A) = \dim_{bc}(\overline{A}).$$

The proof of this proposition is found in [11], and some application of it can be found in D. H. E. Gross's paper [15]. This proposition is important because this is the key to determine the difference between the Hausdorff dimension and the box-counting dimension. This proposition doesn't always hold for the Hausdorff measure. Let's look at an example.

**Example 3.5.5.** *Let $Q = \mathbb{Q} \cap [0,1]$. Given some $\epsilon > 0$, we need $1/\epsilon$ $\epsilon$-covers*

for $Q$ since $\mathbb{Q}$ is dense in $[0, 1]$. Then using equation (3.1),

$$\lim_{\epsilon \to 0} \frac{1/\epsilon}{1/\left(\epsilon^{d_{bc}}\right)} \qquad (3.2)$$

implies if $d_{bc} < 1$, then (3.2) $\to \infty$. If $d_{bc} > 1$, (3.2) $\to 0$. So we obtain the critical point $\dim_{bc}(Q) = 1$.

Now consider $\overline{Q} = [0, 1]$. Given $\epsilon > 0$, we still need $1/\epsilon$ $\epsilon$-covers to cover $\overline{Q}$. Hence we obtain (3.2) as before giving us $\dim_{bc}(\overline{Q}) = 1$.

## 3.6 Measure Properties

First, we discuss general properties of a dimension. Some of these properties will hold for all the dimensions that which we have mentioned, while others will not hold for all dimensions.

**Proposition 3.6.1.** *Let $A$, $B \subseteq \mathbb{R}^n$. Also assume the $A_i$'s and the singleton set $\{a\}$ are subsets of $\mathbb{R}^n$. The following is a list of properties of a dimension.*

1. *If $A \subseteq B$, then $\dim(A) \leq \dim(B)$.*

2. *$\dim(A \cup B) = \max(\dim A, \dim B)$.*

3. *$\dim \left( \bigcup_{i=1}^{\infty} A_i \right) = \sup_i \dim A_i$.*

4. *$\dim f(A) = \dim(A)$ for a "nice" function $f$.*

5. *$\dim(\{a\}) = 0$.*

6. *For a nonempty open $U \subseteq \mathbb{R}^n$, $\dim(U) = n$.*

7. *$\dim(A) = m$ for obvious $m$-dimensional figure $A \subseteq \mathbb{R}^n$ and $m \leq n$.*

The proof for each of the above properties is very elementary. They are found in K. Falconer's book [9]. Here we will provide a brief explanation why the above properties hold true.

The proof for property 1 is obvious, and property 2 follows from 1. The third property can be proved by using induction and properties 1 and 2. As for number 4, what we mean by a "nice" function is that the function is an isometry which includes reflection, rotation, translation, glides and the identity map. In simple terms, the distance between any two distinct points in $A$ is preserved.

The proof for property number 5 is obvious. This property also implies that all finite and countable sets have dimension 0.

To prove property number 6, we use the definition for a set to be open. For any $x \in U$ there is an open ball of $x$ so that the neighborhood is contained in $U$. The dimension of this open neighborhood is $n$. Since the union of these open balls is contained in $U$, the dimension of $U$ is $n$.

Finally if $A$ is an obvious $m$-dimensional object, then the dimension of $A$ is $m$. An obvious $m$-dimensional object includes a straight, connected line in a plane which has dimension 1. A square in $\mathbb{R}^3$ has dimension 2. A cube in a hyperspace has dimension 3.

Each of these properties holds true for the Hausdorff dimension, but they're not always true for the box-counting dimension. The problem arises from the third property [37]. Let's look at an example.

**Example 3.6.2.** *Let $Q = \mathbb{Q} \cap [0, 1]$ as before. We may rewrite $Q$ as $\bigcup_{i \in \mathbb{N}} \{q_i\}$. Then by part 5 from previous proposition, $\dim_H(\{q_i\}) = 0$ for each $i$. And by part 3,*

$$\dim_H \left( \bigcup_{i \in \mathbb{N}} \{q_i\} \right) = \sup_{i \in \mathbb{N}} \dim_H (\{q_i\}) = 0.$$

*On the other hand,*

$$\dim_{bc}(Q) = \dim_{bc}\left( \overline{\mathbb{Q} \cap [0, 1]} \right) = \dim_{bc}([0, 1]) = 1.$$

Hence the box-counting dimension is not infinitely additive.

# Chapter 4

# Examples Using Standard Hausdorff Methods

In this chapter we compute the Hausdorff dimension of some well-known sets.

## 4.1 Finding the Dimension of "Standard" Sets

We first consider some simple shapes in the Euclidean space.

**Example 4.1.1.** *Let $I$ be the unit interval, and let $0 < \delta \leq 1$. Then we need $1/\delta$ of these intervals to cover $I$. Computing for the Hausdorff dimension,*

$$H^s_\delta(I) = \frac{1}{\delta}(\delta)^s = \left(\frac{\delta^s}{\delta}\right).$$

*This implies that our critical point $s = 1$ where $H^s_\delta(I) = 1$.*

*We can easily see that the Hausdorff measure corresponds with both the Lebesgue measure and the standard approach to compute the length of an interval.*

**Example 4.1.2.** *Next, consider the unit square $S$ in the Euclidean plane. We will suppose the four vertices of $S$ are $(0, 0)$, $(1, 0)$, $(1, 1)$ and $(0, 1)$.*

*Similar as before, let $0 < \delta \leq 1$. We let a collection of smaller squares with side-length $\delta$ (so the diameter should be $\delta\sqrt{2}$) cover $S$. Since we need $1/\left(\delta^2\right)$ squares, we compute*

$$H_\delta^s(S) = \frac{1}{\delta^2}\left(\delta^s\right) = \frac{\delta^s}{\delta^2}$$

*which shows that $s = 2$ and $H_\delta^s(S) = 1$.*

*Now computing for the Lebesgue measure of $S$, we have*

$$\lambda(S) = |(1,\ 0) - (0,\ 0)| \times |(0,\ 1) - (0,\ 0)|$$
$$= |(1,\ 0)| \times |(0,\ 1)| = 1.$$

*where $|(x,\ y)|$ is defined to be $(x^2 + y^2)^{\frac{1}{2}}$. Also, the standard way to find the area of a square is to multiply the length of any two perpendicular sides together.*

We note that if we were to choose a collection of circles with diameter $0 < \delta \leq 1$ to cover $S$, then the above will work similarly. The only change we will notice is that $H_\delta^s(S) = 4/\pi$.

We observe that we can compute the length, area or the volume of simple geometric shapes using Lebesgue or Hausdorff measure, or use standard elementary methods and obtain the same results. But when we begin to analyze fractals, Lebesgue measure tells us nothing (because we obtain Lebesgue measure 0 for fractals), and there is no simple elementary method to compute the dimension or the measure of a fractal. So this is where Hausdorff measure becomes important since it exists for all sets.

## 4.2 Fractals and their Hausdorff Dimension

Hence we begin analyzing fractals and their Hausdorff dimension in this section. The fractals we are about to analyze are well-behaving sets which are classified as self-similar. Since we do not need such details in this section, more details about self-similarity will be discussed in section 6.4.

**Example 4.2.1 (Cantor Set).** *Let $C$ be the Cantor set that satisfies* $\sup\{x, y : x, y \in C\} = 1$. *This set is generated by cutting the middle $1/3$ from each of the previous connected lines. So let $C_0$ be the initial unit interval. Then for each $k^{th}$ iteration we obtain $C_k$. Hence we observe*

$$C_0 = [0, 1]$$
$$C_1 = [0, 1/3] \cup [2/3, 1]$$
$$C_2 = [0, 1/9] \cup [2/9, 3/9] \cup [6/9, 7/9] \cup [8/9, 1]$$
$$\vdots \qquad \qquad \vdots$$
$$C_m = \bigcap_{i=0}^{m} C_i.$$

*At $C_0$, we need one interval of unit length to cover $C$. After its first iteration, we need $2$ intervals of length $1/3$. At its second iteration we need $2^2$ intervals of length $1/(3^2)$. Take this to some $m^{th}$ iteration and we need $2^m$ intervals with length $1/(3^m)$ to cover the $C_m$. So we compute*

$$H^s_\delta(C) = 2^m \left(\frac{1}{3}\right)^{ms}.$$

*which implies our critical point $s$ is $\log 2 / \log 3$.*

**Remark 4.2.2.** *The first-iteration of a self-similar fractal gives enough information to determine the fractal's behavior and its Hausdorff dimension.*

Now we work with self-similar sets whose topological dimension is of at least 1. Because it's a bit complicating describing higher dimensional fractals using intervals, we refer the reader to various texts for a pictorical example. But we will describe its changing perimeter, area and volume (whenever they exist), and then compute its Hausdorff dimension.

**Example 4.2.3 (Koch Curve).** *We refer to [28] for pictorial diagrams and details on the Koch curve $K$.*

*We let $P_0$ be the perimeter of the Koch curve at its initial iteration, $P_1$ is at the first iteration and so forth. We shall obtain the perimeter of the*

Koch curve at some $m^{th}$ iteration so we will have a general idea on the rate of the changing perimeter. So computing the $P_i$'s

$$k = 0 \qquad P_0 = 3 \times 1$$
$$k = 1 \qquad P_1 = P_0 - 3 \times 4^0 \left(\frac{1}{3} \times 1\right) + 3 \times 4^0 \times 2 \left(\frac{1}{3} \times 1\right)$$
$$\vdots \qquad \vdots$$
$$k = m \qquad P_m = 3 \times \left[1 + \frac{1}{2} \sum_{k=1}^{m} \left(\frac{4}{3}\right)^k\right].$$

Because we can also compute the changes on the area of the Koch curve, we denote $A_0$ as the area of the initial iteration, $A_1$ is its first iteration, etc. Computing each $A_i$ we get

$$k = 0 \qquad A_0 = \frac{1}{2}\frac{\sqrt{3}}{2}$$
$$k = 1 \qquad A_1 = A_0 + 3 \times 4^0 \frac{1}{2} \left(\frac{1}{3}\right) \left(\frac{\sqrt{3}}{2}\frac{1}{3}\right)$$
$$\vdots \qquad \vdots$$
$$k = m \qquad A_m = \frac{1}{2}\frac{\sqrt{3}}{2} \left[1 + \frac{3}{4} \sum_{k=1}^{m} \left(\frac{4}{3^2}\right)^k\right].$$

Now we make a covering of $K$. Choose one of the sides of length 1 when $k = 0$ to compute $\dim_H(K)$. At its initial iteration, nothing happens to our set. At iteration 1, the side of length 1 breaks into $1/3$ and we need $4^1$ intervals of length $1/(3^1)$ to cover our edge.

When $k = 2$, we need $4^2$ straight lines of length $1/(3^2)$ to cover the same edge at its second iteration. So at $m^{th}$ iteration, we need $4^m$ intervals of length $1/(3^m)$ to cover our set. So we compute

$$H_\delta^s(K) = 4^m \left(\frac{1}{3}\right)^{ms}$$

where the critical point $s$ equals $\log 4/\log 3$.

Note that we're not covering the interior of the Koch curve, but instead, we're covering the boundary of the fractal. We could have covered this set using triangles but we find covering the boundary much simpler for this case. But if we wanted to compute the Hausdorff dimension by covering the interior of the Koch curve, then we use triangles since they are the optimal covering for $K$.

**Example 4.2.4 (Quadratic Koch Curve).** *The quadratic Koch curve $K_Q$ is very similar to the "standard" Koch curve. First, we refer to H. Lauwerier's text [21] for a more detailed description of $K_Q$. Here, we will compute the perimeter $P_k$ and the area $A_k$ of $K_Q$ at various $k^{th}$ iteration to give a better insight on this fractal's behavior.*

Let $P_k$ be the perimeter of $K_Q$, and let $A_k$ be the area of $K_Q$ at iteration $k$. We will let $P_0$ equal the boundary of a unit square, and let $A_0$ be its corresponding area. The perimeter of the quadratic Koch curve is

$$k = 0 \qquad P_0 = 4$$
$$k = 1 \qquad P_1 = P_0 - 8^1\left(\frac{1}{4}\right) + 8^0 \times 4\left(6 \times \frac{1}{4}\right)$$
$$\vdots \qquad \vdots$$
$$k = m \qquad P_m = 4\left[1 + \sum_{k=0}^{m-1} 2^k\right],$$

while the area of $K_Q$ is

$$k = 0 \qquad A_0 = (1)^2$$
$$k = 1 \qquad A_1 = A_0 - 4 \times 8^0 \left(\frac{1}{4^1}\right)^2 + 4 \times 8^0 \left(\frac{1}{4^1}\right)^2 = A_0 = (1)^2$$
$$\vdots \qquad \vdots$$
$$k = m \qquad A_m = (1)^2.$$

Now computing for the Hausdorff dimension of $K_Q$, we use $16^k$ squares with side-length $\delta = 1/4^k$ to cover the quadratic Koch curve. Hence we obtain

$$H^s_\delta(K_Q) = \left(16^k\right)\left(\frac{1}{4^k}\right)^s$$

and we find $s = \log 16 / \log 4 = 2$.

We notice that the area of $K_Q$ is not affected by the changes of its boundary, but we will continue to refer to the quadratic Koch curve as a fractal [21].

**Conjecture 4.2.5.** *The Hausdorff dimension of the quadratic Koch curve may have a direct relationship to its area. Since the area of $K_Q$ is an integer, the Hausdorff measure may not recognize the way in which the points of $K_Q$ are scattered on the Euclidean plane.*

Further analysis on the relationship between fractals and their Hausdorff dimension is recommended. K. J. Falconer classifies sets under three different categories: one for sets with Hausdorff dimension $s$ that's strictly between 0 and 1, second for sets with $s > 1$ and the third section for sets in higher dimensions. He discusses that any set with Hausdorff dimension less than 1 is irregular, totally disconnected, the density fails to exist at almost every point of the set, etc. We refer the text "The Geometry of Fractal Sets" [8] for more details.

**Example 4.2.6 (Sierpinski Gasket).** *The Sierpinski gasket $S_G$ is another name for the Sierpinski triangle. We take the reader to [1] or [2] for instructions on generating this fractal. Similar as before, we compute $P_k$ and $A_k$. Let $P_0$ be the initial perimeter of length 3 and $A_0$ be the area $(\sqrt{3})/4$. Then the perimeter changes by*

$$
\begin{aligned}
k &= 0 & P_0 &= 3 \\
k &= 1 & P_1 &= P_0 + 3^0 \left(3 \times \frac{1}{2}\right) \\
&\vdots & &\vdots \\
k &= m & P_m &= 3\left[1 + \frac{1}{2}\sum_{k=0}^{m-1}\left(\frac{3}{2}\right)^k\right],
\end{aligned}
$$

and the area computes as

$$k = 0 \qquad A_0 = \frac{1}{2}\frac{\sqrt{3}}{2}$$

$$k = 1 \qquad A_1 = A_0 - 3^0 \frac{1}{2}\left(\frac{1}{2^1}\right)\left(\frac{\sqrt{3}}{2}\frac{1}{2}\right)$$

$$\vdots \qquad \vdots$$

$$k = m \qquad A_m = \frac{1}{2}\frac{\sqrt{3}}{2}\left[1 - \frac{1}{3}\sum_{k=1}^{m}\left(\frac{3}{2^2}\right)^k\right].$$

We observe that $3^m$ triangles of side-length $1/(2^m)$ are needed to cover the Sierpinski gasket at each $m^{th}$ iteration. So

$$H_\delta^s(S_G) = 3^m \left(\frac{1}{2^m}\right)^s$$

which implies $s = \log 3/\log 2$.

**Example 4.2.7 (Sierpinski Carpet).** We will refer to the websites [3] and [36] for some diagrams of the Sierpinski carpet $S_C$, but we will briefly mention the construction of $S_C$. We begin with a unit square. Divide the unit square into 9 squares of equal size. Each square has side-length $1/3$. Delete the middle square. This was the first iteration. The second iteration repeats this process on each of our smaller squares which should require us to cut-out 8 squares of side-length $1/(3^2)$. Continue with this iteration to obtain the Sierpinski carpet.

Let's compute the perimeter $P_k$ and the area $A_k$ of $S_C$ at various iterations. Let $P_0$ and $A_0$ be the initial perimeter and area of $S_C$ respectively. Then

$$k = 0 \qquad P_0 = 4$$

$$k = 1 \qquad P_1 = P_0 + 8^0 \times 4\left(\frac{1}{3}\right)$$

$$\vdots \qquad \vdots$$

$$k = m \qquad P_m = 4\left[1 + \frac{1}{8}\sum_{k=1}^{m}\left(\frac{8}{3}\right)^k\right]$$

and

$$k = 0 \qquad A_0 = (1)^2$$
$$k = 1 \qquad A_1 = A_0 - 8^0 \left(\frac{1}{3}\right)^2$$
$$\vdots \qquad \vdots$$
$$k = m \qquad A_m = (1)^2 \left[1 - \frac{1}{8} \sum_{k=1}^{m} \left(\frac{8}{3^2}\right)^k\right].$$

Hence for each $m$, we need $8^m$ squares of side-length $1/(3^m)$ to cover the Sierpinski carpet. So

$$H_\delta^s (S_C) = 8^k \left(\frac{1}{3^k}\right)^s,$$

and $s = \log 8 / \log 3$.

**Example 4.2.8 (Menger Sponge).** *Let's denote $M_S$ as the Menger sponge. A diagram and a description is found in [35]. We will compute the surface area $A_k$ and the volume $V_k$ of $M_S$ at different iterations. Let $A_0$ equal 6 since we have 6 surfaces of area that which is of a unit square, and we let $V_0$ be $(1)^3$ where $V_0$ is the volume of the unit cube before any iterations. So we obtain the following areas*

$$k = 0 \qquad A_0 = 6$$
$$k = 1 \qquad A_1 = A_0 - 6 \times 8^0 \left(\frac{1}{3}\right)^2 + 20^0 \times 6 \times 4 \left(\frac{1}{3}\right)^2$$
$$\vdots \qquad \vdots$$
$$k = m \qquad A_m = 6 \left[1 + \sum_{k=1}^{m} \left(\frac{1}{5}\left(\frac{20}{3^2}\right)^k - \frac{1}{8}\left(\frac{8}{3^2}\right)^k\right)\right]$$

*and the volumes*

$$k = 0 \qquad V_0 = (1)^3$$

$$k = 1 \qquad V_1 = V_0 - 20^0 \left(\frac{1}{3}\right)^3$$

$$\vdots \qquad \qquad \vdots$$

$$k = m \qquad V_m = (1)^3 \left[1 - \frac{1}{20} \sum_{k=1}^{m} \left(\frac{20}{3^3}\right)^k \right].$$

*Since we need $20^m$ sets with diameter $\delta = 1/3^m$ to cover the Menger sponge,*

$$H^s_\delta(M_S) = (20^m) \left(\frac{1}{3^m}\right)^s$$

*and the Hausdorff dimension of $M_S$ is $\log 20 / \log 3$.*

Now having seen some standard fractals let's analyze a hyperfractal.

**Remark 4.2.9 (To Generate a Hypergasket).** *A hypergasket can be generated in the following way:*

1. *take any $n - 1$ dimensional set;*

2. *bring the midpoint of the set, relatively equidistant from every other element in the set, into the $n^{th}$ dimension;*

3. *create additional edges by taking all the vertices in $n - 1$ dimension and extending them to the newly positioned midpoint;*

4. *make faces by taking all edges from $n-1$ dimensional set and extending them to this midoint.*

If we start out with a point we can extend it to another point so that it becomes a line. Then we may extend the midpoint of this line in a perpendicular direction to make an equilateral triangle. Similarly take the midpoint of the equilateral triangle and pull it to its third dimension to obtain a tetrahedron. So we consider a hyperpyramid.

**Example 4.2.10 (Hyperpyramid).** *We refer the reader to C. T. Ho and S. L. Johnsson's paper [16] for a description of hyperpyramids. But one should be able to imagine having a 3-dimensional cube and a point which is relatively equidistant from all the other points, and then draw edges and surfaces as described above. We shall denote the hyperpyramid as $H_P$.*

*We define its first iteration by dividing each edge by 2 and then filling up this set with smaller hyperpyramids. The second iteration is exactly the same as the first iteration except that we repeat this process on each of the smaller hyperpyramids. We find that we need $9^k$ sets of diameter $1/(2^k)$ to cover our $H_P$. Hence*

$$H_\delta^s(H_P) = 9^k \left(\frac{1}{2^k}\right)^s$$

*implies $s = \log 9 / \log 2$.*

# Chapter 5

# Nonstandard Measure Theory

We develop the notion of computing measure in the nonstandard universe. So in this chapter we will explore various versions to measure a set.

From this point we will fix our space $X$ to be the unit interval $[0, 1]$ to simplify our analysis unless we explicitly say otherwise. We let $\Omega$ be

$$\left\{0, \ldots, \frac{i}{N}, \ldots, \frac{N}{N}\right\}$$

where $N$ is a nonstandard natural. As we can see, the standard map of $\Omega$ takes us back to $[0, 1]$.

We recall from the chapter on infinitesimals that $\operatorname{st}^{-1}(A) \subseteq \Omega$ is not in general internal for some set $A \subseteq [0, 1]$. So in particular, $\operatorname{card}\left(\operatorname{st}^{-1}(A)\right)$ is meaningless. Hence we will approximate this set by internal sets.

## 5.1 Nonstandard Version of the Lebesgue Measure

We first discuss how we may go about in obtaining the nonstandard version of the Lebesgue measure. So let $A \subseteq [0, 1]$ be a nonempty set, and let $B$ be an internal subset of $\Omega$ so that $B \subseteq \mathrm{st}^{-1}(A)$. We note that there may be many internal $B$'s contained in $\mathrm{st}^{-1}(A)$, and this collection includes the empty set.

Each $B$ has a discrete measure in $\Omega$ that which is the discrete probability measure on $\Omega$ where each $x \in \Omega$ is equally likely. Let's define the discrete measure of $B$ as

$$d(B) = \frac{\mathrm{card}(B)}{N+1}$$

where $N \in {}^*\mathbb{N}$. This is in ${}^*\mathbb{Q} \cap {}^*[0, 1]$.

**Definition 5.1.1.** *Let $A$ be a subset of $X = [0, 1]$ and let $\Omega$ be the set $\left\{0, \ldots, \frac{i}{N}, \ldots, \frac{N}{N}\right\}$ where $N \in {}^*\mathbb{N} \setminus \mathbb{N}$. Let $d(B)$ be defined as above for some internal set $B$. Then the lower Lebesgue measure of $A$ is*

$$\underline{\lambda}(A) = \sup\{\mathrm{st}\, d(B) : \text{internal } B \subseteq \mathrm{st}^{-1}(A)\},$$

*and the upper Lebesgue measure of $A$ is*

$$\overline{\lambda}(A) = \inf\{\mathrm{st}\, d(B) : \text{internal } B \supseteq \mathrm{st}^{-1}(A)\}.$$

*We say $A$ is Lebesgue measurable if $\overline{\lambda}(A) = \underline{\lambda}(A)$, and we write $\lambda(A)$ for the common value.*

Strictly speaking, this definition is due to Loeb who first showed it is equivalent to Lebesgue's definition. See Cutland's book [6] for details. Further studies on Lebesgue measure may be found in K. J. Falconer's text [8] but the above will suffice for this report.

## 5.2 Discrete Hausdorff $s$-dimensional Outer Measure

Now a nonstandard approach to Hausdorff measure is discussed along similar lines.

First we will mention a nice property of Hausdorff measure. The Hausdorff outer measure is related to the Lebesgue outer measure by

$$H^s(A) = c_s \bar{\lambda}(A)$$

for some constant $c_s$ that depends on $s$, and for nice sets $A$ [6]. They are related to each other in the sense that given a nice set with positive integral dimension $s$ embedded in $\mathbb{R}^n$, we can take this set to the $\mathbb{R}^s$-space and use Lebesgue outer definition to find the measure of $A$. Because of this we may hope to obtain a nonstandard presentation of Hausdorff measure analogous to that for Lebesgue measure.

We now give the nonstandard discrete version of Hausdorff measure that we require.

**Definition 5.2.1.** *Given $\delta > 0$ and $s \in [0,1]$, a $\delta$-interval of $\Omega$ is a set $\{i/N, (i+1)/N, \ldots, j/N\}$ with diameter $(j-i+1)/N \leq \delta$. For an internal set $B \subseteq \Omega$ let the discrete $s$-dimensional measure be defined as*

$$h_\delta^s(B) = \min \sum_{i=1}^{L} (\mathrm{diam}(V_i))^s$$

*where we take the minimum over all partitions $\{V_1, \ldots, V_L\}$ of $B$ into $\delta$-intervals.*

Here our nonstandard finite $L \in {}^*\mathbb{N}$ varies according to each partition. Note that in the discrete space $\Omega$, this definition is internal and hence (because there are nonstandard finitely many partitions to consider) the min is a true minimum. Also since our $V_i$'s are subsets of $\Omega$ which is already normalized, we do not need to normalize it again.

Next we say how this discrete measure gives rise to Hausdorff $s$-dimensional outer measure.

**Theorem 5.2.2.** *Let $A$ be a subset in $[0,1]$ and let $H^s(A)$ be the Hausdorff $s$-dimensional outer measure. Then*

(a) $H^s(A) \leq \inf\{\operatorname{st} h_\delta^s(B) : \text{internal } B \supseteq \operatorname{st}^{-1}(A)\}$ *for all infinitesimal $\delta$*

(b) $H^s(A) = \lim_{\delta \to 0} \inf\{\operatorname{st} h_\delta^s(B) : \text{internal } B \supseteq \operatorname{st}^{-1}(A)\}$ *where, in the limit, $\delta$ ranges over standard positive values from $\mathbb{R}$.*

In particular the nonstandard version of $H^s(A)$ does not depend on $N$.

The following proof has been obtained after some discussions, advice and suggestions from my supervisor Dr. Richard W. Kaye.

*Proof.* For (a), we assume that $\delta$ is infinitesimal and fixed, and $h_\delta^s(B)$ is finite for some internal $B \supseteq \operatorname{st}^{-1}(A)$. Let $\eta > 0$ be standard. We must find an $\eta$-cover $U_i$ for $i \in \mathbb{N}$ such that $\sum_i (\operatorname{diam}(U_i))^s \leq \operatorname{st}(h_\delta^s(B))$.

Take an optimal partition $\{V_j\}_{j=1}^K$ of $B$ in the sense of $h_\delta^s(B)$. Now, using an internal induction in the nonstandard universe we modify $\{V_j\}_{j=1}^K$ to some other partition $\{W_k\}_{k=1}^L$ as follows. When defining some $W_k$, given an internal interval $I \subseteq \Omega$ consider the leftmost $V_i$ that is to the right of $I$. If $I \cup V_i$ is also an interval, we replace $I$ with this interval and carry on. This process stops when one of two things happens:

- $I \cup V_i$ is not an interval, or there is no further $V_i$ to the right of $I$;

- $\operatorname{diam}(I) \geq \eta$.

If either of these happens we stop the construction and let $W_k = I$. We then construct $W_{k+1}$ starting with $I'$, the leftmost $V_i$ to the right of $I$ if any.

At the end of this construction we will have nonstandard finitely-many intervals $W_k$ of length at most $\eta + \delta$ partitioning $B$. This gives a countable $\eta$-cover $\mathcal{U}$ of our original $A$ consisting of all sets $U = \operatorname{st}(W_k)$ for some $W_k$, such that the set $U$ has nonempty interior.

To see that this does indeed cover $A$, consider $a \in A$. Then the monad of $a$ is contained in $B$. By overspill, this monad is contained in an interval $J \subseteq B$ where the length of $J$ is not infinitesimal. But then, by construction, the monad of $a$ is contained in either some set $W_k$ or else in two neighboring sets $W_k$ and $W_{k+1}$. (This latter case is when the construction of the set $W_k$ was 'finished' whilst inside the monad of $a$.) Moreover $W_k$ (or in the other case both of $W_k$ and $W_{k+1}$) have lengths that are non-infinitesimal, by construction and by choice of $\eta$. So $a$ is either in the interior of $\text{st}(W_k)$ or else is an endpoint of both $\text{st}(W_k)$, $\text{st}(W_{k+1})$, as required. The cover $\mathcal{U}$ is countable because any set of intervals in $\mathbb{R}$ all with nonempty interior must necessarily be countable. Finally, since without loss $s \leq 1$ (since we are working on the unit interval) and hence $(a+b)^s \leq a^s + b^s$ for $a, b > 0$, we have

$$\sum_{U \in \mathcal{U}} (\text{diam } U)^s = \text{st} \sum_k (\text{diam}(W_k))^s \leq \text{st} \sum_i (\text{diam } V_i)^s,$$

as required.

For (b), we observe first that for $0 < \delta < \eta$, even if $\delta$ is not infinitesimal, then the argument just given shows that $H^s_\eta(A) \leq \text{st } h^s_\delta(B)$ for all internal $B \supseteq \text{st}^{-1}(A)$. Thus we only have to prove the other direction. It suffices to show that for each standard $\delta > 0$ and each standard $\epsilon > 0$ there is an internal $B \supseteq \text{st}^{-1}(A)$ with $H^s_\delta(A) \geq \text{st}(h^s_\eta(B)) - \epsilon$ for a certain $\eta$ depending on $\delta$ and $\epsilon$ only. Let $\lambda = H^s_\delta(A)$.

Let $\mathcal{U} = \{U_i\}_{i \in \mathbb{N}}$ be a $\delta$-cover of $A$ so that

$$\sum_{i=1}^N (\text{diam}(U_i))^s \leq \lambda + (\epsilon/2),$$

and assume without loss that each $U_i$ is an interval. We 'enlarge' $U_i$ by increasing its length by $(\epsilon/2)^{1/s} 2^{-i/s}$ on each side, obtaining intervals $V_i$. By saturation there is a sequence of intervals $W_i \subseteq \Omega$ such that $\text{st}(W_i) = V_i$ and $W_i$ is defined for all $i < K$, where $K > \mathbb{N}$. Then for any $N > \mathbb{N}$ we have $\bigcup_{i=1}^N W_i \supseteq \text{st}^{-1}(A)$ because of the 'enlarging'. Moreover for each $N \in \mathbb{N}$ we

have

$$\sum_{i=1}^{N}(\operatorname{diam}(W_i))^s \leq \sum_{i=1}^{N}(\operatorname{diam}(U_i) + (\epsilon/2)^{1/s}2^{-i/s})^s$$
$$\leq \sum_{i=1}^{N}(\operatorname{diam}(U_i))^s + (\epsilon/2) + (\epsilon/2)\sum_{i=1}^{N}2^{-i}$$
$$\leq \lambda + \epsilon.$$

So by overspill there is an infinite $N$ such that this holds. Thus some $B = \bigcup_{i=1}^{N} W_i$ has a partition showing $h_\eta^s(B) \leq \lambda + \epsilon$. Finally note that the maximum diameter of any $W_j$ is $\delta + (\epsilon/2)^{1/s}2^{-1/s}$ which may be made as close to $\delta$ as we like by choosing $\epsilon$ sufficiently small. This completes the proof. □

Unfortunately it does not seem possible to replace the limit as $\delta \to 0$ over standard $\delta$ with an infinitesimal $\delta$ in Theorem 5.2.2. In other words, the inequality in part (a) of the theorem cannot be replaced by an equality, even for carefully chosen $\delta$, as the example of the Cantor set shows.

**Example 5.2.3.** *Consider the Cantor set $C \subseteq [0,1]$ which has dimension $s = \log 2/\log 3$. If $B \supseteq \operatorname{st}^{-1}(C)$ is internal, and $\delta > 0$ is infinitesimal, we try to estimate $h_\delta^s(B)$, giving a lower bound.*

*For each $a \in C$ the monad $\operatorname{st}^{-1}(a)$ is covered by an interval of non-infinitesimal length $I \subseteq B$. Taking all such intervals and mapping back to $[0,1]$ via the standard part map, and taking the interiors of these sets, we obtain an open cover of $A$. Since $C$ is bounded and closed hence compact, there is a finite subcover, and this shows that without loss we may assume that $B$ is one of the sets of $2^m$ intervals of length $3^{-m}$ obtained at the $m^{th}$ stage of the construction of $C$. (If not some such collection would be smaller than $B$.) For an interval $I \subseteq \Omega$ of length $\ell$ we have that $h_\delta^s(I)$ is approximately $(\ell/\delta)\delta^s$, so*

$$h_\delta^s(B) = 2^m \frac{3^{-m}}{\delta}\delta^s = (2/3)^m \delta^{s-1}.$$

*But because $0 < s < 1$ and $\delta$ is infinitesimal, this value is infinite for all standard $m \in \mathbb{N}$. Hence, $\operatorname{st}(h_\delta^s(B)) = \infty$ for all internal $B \supseteq \operatorname{st}^{-1}(C)$.*

We can speculate that the problem here is that for infinitesimal $\delta$, the function $h_\delta^s$ measures the size of $B$ in terms of $s < 1$, whereas the local structure of $B$ is that has dimension one. It would be nice to be able to investigate this more, asking perhaps for which sets $A$ is it true that the inequality in (a) of Theorem 5.2.2 becomes an equality for suitably chosen infinitesimal $\delta$.

## 5.3 Nonstandard Box-Counting Conjecture

In this section we shall employ nonstandard techniques to obtain the box-counting dimension of a set.

Let $A$ be a subset of $[0, 1]^n$. Suppose $0 < h \ll k \approx 0$ so that $mh < k$ for all standard $m \in \mathbb{N}$. We now subdivide $^*[0, 1]^n$ into neighborhoods by defining a neighborhood of hypervolume $h^n$ by

$$B_{\underline{x}}^h = \{\underline{y} \in {}^*[0, 1]^n : x_i \leq y_i < x_i + h \text{ for all } i\} \qquad (5.1)$$

around each $\underline{x} \in {}^*[0, 1]^n$. Similarly we define a neighborhood of hypervolume $h^n$ around each $\underline{x} \in {}^*[0, 1]^n$ by setting

$$B_{\underline{x}}^k = \{\underline{y} \in {}^*[0, 1]^n : x_i \leq y_i < x_i + k \text{ for all } i\}. \qquad (5.2)$$

These neighborhoods will be used to 'sample' our set $A$ at intervals of $h$ and $k$. More precisely, let $A_h$ and $A_k$ be the sets defined by

$$A_h = \left\{\underline{x} = (i_1 h, i_2 h, \ldots i_n h) : \underline{i} \in {}^*\mathbb{N},\ \underline{x} \in {}^*[0, 1]^n \text{ and } B_{\underline{x}}^h \cap {}^*A \neq \varnothing\right\} \qquad (5.3)$$

and

$$A_k = \left\{\underline{x} = (i_1 k, i_2 k, \ldots i_n k) : \underline{i} \in {}^*\mathbb{N},\ \underline{x} \in {}^*[0, 1]^n \text{ and } B_{\underline{x}}^k \cap {}^*A \neq \varnothing\right\}, \qquad (5.4)$$

respectively. Let $|A_h|$ and $|A_k|$ denote the number of elements in $A_h$ and $A_k$ respectively. We are interested in the growth rate of $|A_h|$ as $h$ becomes small. Put another way, we are interested in the comparative sizes of $|A_h|$ and $|A_k|$ for $0 < h \ll k \approx 0$. So we define

$$d_{h,\,k}(A) = \operatorname{st}\left(\frac{\log\left(\frac{|A_h|}{|A_k|}\right)}{\log\left(\frac{k}{h}\right)}\right). \tag{5.5}$$

With these definitions we have

**Conjecture 5.3.1 (Nonstandard Box-Counting Measure).** *Let $A$ be a subset of $[0,\,1]^n$, and let $0 < h \ll k \approx 0$ be infinitesimal. With the notation given above, if $d_{bc}(A) = d_{bc}$ is the classical box-counting dimension of $A$ then $d_{bc}(A) = d_{h,\,k}(A)$. In particular, $d_{h,\,k}(A)$ does not depend on the choice of $h, k$.*

When $h$, $k$, $A$ are understood from the context we shall abbreviate $d_{h,\,k}(A)$ as $d$.

Now, we make some remarks that hopefully should help motivate the theorem and its proof.

Firstly, keep in mind that the grids of interval size $h$ and $k$ do not have to match the holes or the discontinuities in our set $A$, for example, if it is one of the standard fractal examples defined by using a regular and decreasing grid size. So the sampling process only promises there is a point in $^*A$ nearby: in other words even though $B_{\underline{x}}^h$ might not be contained in $^*A$, as long as some $\underline{x} \in {}^*[0,\,1]^n$ is in $^*A$, the cover $B_{\underline{x}}^h$ will contribute to $|A_h|$. It is similar for the neighborhoods $B_{\underline{x}}^k$ and the set $A_k$. Also, if the grid sizes $h$ and $k$ are very close to each other, then it is possible that the ratio of the cardinalities of $A_h$ and $A_k$ may not indicate the required growth rate. Indeed the cardinality of $A_h$ may even be the same as the cardinality of $A_k$ because of some peculiarity to do with the choice of $k$ and $h$. To get something more interesting reflecting the nature of $A$ we need to put some conditions on $h$ and $k$.

Secondly, since the grids of intervals $h$, $k$ are not required to match the set $^*A$, the sets $A_h$ and $A_k$ will contain some $\underline{x}$ that may or may not be in the set $^*A$. But since both $0 < h, k$ are infinitesimals, $|A_h|$ and $|A_k|$ are in $^*\mathbb{N}$, and

$$\frac{\log\left(\frac{|A_h|}{|A_k|}\right)}{\log\left(\frac{k}{h}\right)} \gtrapprox d_{bc},$$

where the left-hand side might be greater than $d_{bc}$ by an infinitesimal amount. But this is okay since when we take the standard part of the left-hand side, this number equals the classical box-counting dimension.

Thirdly, our defined number $d = d_{bc}$ should be independent of the grid size, the imbedding and rotation, translation and reflection of the set $A$. To say that $d$ is independent of $h$ and $k$, we mean that we cannot see the difference between $h$ and $k$. All we know is that they are both some small, if not exact, number in the real world.

Finally, note that the notion of a "sample" of our set $A$ is the intuitive one of selecting points from $A$ according to how this set meets the squares, cubes or hypercubes in our grid. This is the method commonly used in primary school to estimate the area of an irregular shape. Here it is being used to estimate some other property of our shape, but the relationship to measure should be clear. By sampling on an infinitesimal grid, we may lose some information about $A$, but we can at least recover the closure of $A$ since $\text{st}(^*A) = \overline{A}$ which will be proved in Proposition 5.3.2. This makes our conjecture very similar to the box-counting dimension since the box-counting dimension of the closure of $A$ equals $\dim_{bc}(A)$. This is one of the main reasons why our conjecture cannot be the nonstandard version of the Hausdorff measure.

The following proposition uses a grid of an infinitesimal step-size $h$ over a set $A$, and then determines the standard part of $A_h$, which is the closure of $A$. The notation for nonstandard elements and sequences slightly changed in this proposition but it should be fundamentally easy to follow through.

**Proposition 5.3.2.** *Let $h > 0$ be infinitesimal and $A \subseteq [0, 1]^n$. Then with $A_h$ as defined in (5.3), we have $\overline{A} = \{\text{st}(\underline{x}): \underline{x} \in A_h\}$.*

*Proof.* Let $B = \{st(\underline{x}) : \underline{x} \in A_h\}$. If $x \notin \overline{A}$ then there is an open neighborhood $U$ around $x$ so that $U \cap A = \varnothing$. So since $\underline{x} \notin {}^*U \cap {}^*A$, $x$ is not contained in $B$.

If $x \in \overline{A}$ then there is a sequence $S = \{x^i\}_{i \in \mathbb{N}}$ that converges to $x$ where $x^i \in A$ for all but finitely-many $i$'s. So since

$$\forall \, \epsilon > 0 \,\exists\, N > 0 \left(x^i \in A \,\&\, \left|x^i - x\right| < \epsilon \,\forall\, i > N\right)$$

is a first-order sentence, by Łoś' theorem

$$\forall \, \epsilon > 0 \,\exists\, N > 0 \,\&\, \left(\underline{x}^i \in {}^*A \,\&\, \left|\underline{x}^i - \underline{x}\right| < \epsilon \,\forall\, i > N \,\&\, i \in {}^*\mathbb{N}\right)$$

holds true in the nonstandard world where $\underline{x}^i = [x_0^i,\, x_1^i,\, x_2^i, \ldots]$.

So for any open set ${}^*U$ of $\underline{x}$, almost all hyperreal numbers in the sequence ${}^*S = \{\underline{x}^0,\, \underline{x}^1,\, \underline{x}^2 \ldots\}$ are contained in ${}^*U \cap {}^*A$. Hence $x \in B$.

Finally, it is clear that $B \subseteq \overline{A}$ since for any $\underline{x} \in A_h$

$$B = \text{st}\left(B_{\underline{x}}^h \cap {}^*A\right) \subseteq \text{st}({}^*A) = A \subseteq \overline{A}.$$

Hence the closure of $A$ is the set $\{\text{st}(\underline{x}) : \underline{x} \in A_h\}$. □

### 5.3.1 Examples

Let's apply equation (5.5) to some examples.

**Example 5.3.3.** *Let $A = [0, 1]$ and $0 < h \ll k \approx 0$. Let*

$$A_h = \{ih : 0 < ih < 1 \text{ and } i \in {}^*\mathbb{N}\}$$

*and let*

$$A_k = \{ik : 0 < ik < 1 \text{ and } i \in {}^*\mathbb{N}\}.$$

We define the bracket $[x]$ as the integer part of $x$. Then $|A_h| = \left[\frac{1}{h}\right] - 1$ and $|A_k| = \left[\frac{1}{k}\right] - 1$ where we subtract one or two of the endpoints depending on the set $A$. So we compute

$$d_{h,\,k}(A) \approx \text{st}\left(\frac{\log(|k|/|h|)}{\log(|k|/|h|)}\right) = \text{st}(^*1) = 1.$$

The subtraction of the endpoints from $A_h$ and $A_k$ does not significantly affect the sizes of $A_h$ and $A_k$ since they contribute very little compared to $\left[\frac{1}{h}\right]$ and $\left[\frac{1}{k}\right]$. In fact, the constants disappear when we take the standard part of $A_h$ and $A_k$.

**Example 5.3.4.** *Similarly, let $A = [0,\,1]^2$ and we will let $0 < h \ll k \approx 0$. Suppose*

$$A_h = \{(ih,\,jh) : 0 < ih,\,jh < 1 \text{ and } i,\,j \in {}^*\mathbb{N}\}$$

*and*

$$A_k = \{(ik,\,jk) : 0 < ik,\,jk < 1 \text{ and } i,\,j \in {}^*\mathbb{N}\}.$$

*Then the size of $A_h$ is $\left(\left[\frac{1}{h}\right]^2 - a\right)$ and the size of $A_k$ is $\left(\left[\frac{1}{k}\right]^2 - b\right)$ where $a$, $b \in \mathbb{N}$. Again since $a$ and $b$ are small compared to $\left[\frac{1}{h}\right]^2$ and $\left[\frac{1}{k}\right]^2$, they will disappear when we take the standard part. Hence*

$$d_{h,\,k}(A) \approx \text{st}\left(\frac{\log(|k|^2/|h|^2)}{\log(|k|/|h|)}\right) = \text{st}(^*2) = 2.$$

**Example 5.3.5 (Sierpinski Carpet).** *Given $S_C$ as before, let's put two grids of step-sizes $0 < h \ll k \approx 0$ over our set. Let's put the following conditions on this example: the sizes of $h$ and $k$ are chosen so that the fraction $k/h$ is infinite; the grid with step-size $h$ is lined up in the same direction as that of the $k$-grid; and the grid with step-size $k$ is lined up so that each grid area of side-length $k$ contains an exact but smaller copy of $S_C$.*

*Let's denote $S_{Ch}$ as $S_C$ with the $h$-grid over the set, and let $S_{Ck}$ be $S_C$ with the $k$-grid over our set. We zoom-in on one of the boxes of $S_{Ck}$ to* calculate d. We cut out $(8^{N-1})$-squares at each iteration, but this time, each square has side-length $\left(k \times \frac{1}{3}^N\right)$ where $N$ ranges over $0 \leq N \leq \left[\log_3\left(\frac{k}{h}\right)\right]$.

We have taken $\log_3$ since the side-length is reducing by a factor of $\frac{1}{3}$. Also we shall take the integer part of $\log_3\left(\frac{k}{h}\right)$ if necessary. Now we do some computation.

We obtain $|S_{Ck}| = 1$ while

$$|S_{Ch}| = \left(\frac{k}{h}\right)^2 \times \left(1 - \sum_{i=1}^{N} \frac{8^{i-1}}{9^i}\right) = \left(\frac{k}{h}\right)^2 \times \left(\frac{8}{9}\right)^N,$$

where the fraction $k/h$ represents the relative scale factor, and $(8/9)^N$ is the fraction of the area that remains at each $N^{th}$ iteration.

So choosing $N$ to be as big as possible, which is when $N = \left[\log_3\left(\frac{k}{h}\right)\right]$, we compute

$$\frac{|S_{Ch}|}{|S_{Ck}|} \approx \left(\frac{8}{9}\right)^N \times \left(\frac{k}{h}\right)^2$$
$$= \left(\frac{k}{h}\right)^{2+N \times \frac{\log(8/9)}{\log(k/h)}}$$
$$= \left(\frac{k}{h}\right)^{2+\frac{\log_3(k/h) \times \log_3(8/9)}{\log_3(k/h)}}$$
$$= \left(\frac{k}{h}\right)^{2+\log_3(8/9)} = \left(\frac{k}{h}\right)^{\log_3 8}.$$

Hence computing for $d_{h,\,k}$,

$$d_{h,\,k}(A) = \text{st}\left[\left(\log \frac{|S_{Ch}|}{|S_{Ck}|}\right) \Big/ \log\left(\frac{k}{h}\right)\right]$$
$$= \text{st}\left[\left(\log \left(\frac{k}{h}\right)^{\log_3 8}\right) \Big/ \log\left(\frac{k}{h}\right)\right]$$
$$= \log_3 8 = \frac{\log 8}{\log 3}.$$

This is the result that which we wanted.

It appears that the analyzed set $A$ must be self-similar when using this conjecture, or else, we must analyze each nonempty $B_{\underline{x}}^k \cap^* A$ because of its random distribution over $\mathbb{R}^n$. This bring us to the simplified nonstandard box-counting approach which doesn't necessary restrict our sets to be self-similar.

Self-similarity is discussed in section 6.4.

## 5.4 Simplified Nonstandard Box-Counting Conjecture

In this section we simplify the 2-grid conjecture and propose the following.

**Conjecture 5.4.1 (Simplified Discrete Box-Counting Conjecture).**
*Let $A$ be a nonempty subset of $[0, 1]^n$, and let $d_{bc}$ be its standard box-counting dimension of $A$. We put a grid of step-size $0 < h \approx 0$ over $A$, and define $A_h$ similar as before. Then*

$$\operatorname{st}\left(\liminf_{h \to 0} |A_h|\, h^d\right) = \begin{cases} \infty & \text{if } d < d_{bc}, \\ 0 & \text{if } d > d_{bc}. \end{cases} \tag{5.6}$$

**Example 5.4.2.** *Let $h \approx 0$ be positive. We let $B_{\underline{x}}^h$ be (5.1) and $A_h \in {}^*\mathbb{N}$ be the set as defined by (5.3). If we choose $A$ to be an interval $(a, b) \subseteq [0, 1]$ then we have*

$$|A_h| = \mathscr{O}\left(\frac{1}{h}\right)$$

*where $\mathscr{O}\left(\frac{1}{h}\right)$ implies that the size of $A_h$ is in some order of $1/h$. That is, $|A_h| = (C \times (1/h) + a)$ where $C$ is a constant and $a$ is an infinitesimal.*

*For a nonempty subset $A$ of $[0, 1]^2$,*

$$|A_h| = \mathscr{O}\left(\frac{1}{h^2}\right),$$

*or we write $|A_h| = (C \times (1/h^2) + a)$ for some constant $C$ and $0 < a \approx 0$.*

*In general if $A$ is a subset of $[0, 1]^d$, then*

$$|A_h| = \mathscr{O}\left(\frac{1}{h^d}\right).$$

### 5.4.1 Examples

The idea behind the simplified box-counting conjecture is that we may replace each relatively connected line, surface or a volume with a dot and we simply "count" the dots to find the dimension each set.

All fractals are closed and the fractals we will examine in this section are self-similar, so the box-counting dimension exists and equals the Hausdorff dimension [21] and [33]. Here we will explore various fractals and find its dimension using our simplified box-counting methods (Conjecture 5.4.1).

**Remark 5.4.3.** *The terminology "total points" or "possible dots" used in this section reflects the number of (connected) subsets which depends on iteration k that which we need to cover the convex closure of A, which is normally $A_0$.*

**Example 5.4.4 (Cantor Set).** *Let's return to the Cantor set $C$ where $\sup\{x, y\} = 1$ for any $x, y \in C$. For each iteration $k$, we shall replace each connected line with a dot to obtain*

$$
\begin{array}{cc}
k = 0 & . \\
k = 1 & .\qquad . \\
k = 2 & .\ \ .\quad .\ \ . \\
k = 3 & .\,.\,.\,.\ \ .\,.\,.\,. \\
\vdots & \vdots
\end{array}
$$

*So we observe that we have $2^k$ points out of a total of $3^k$ for each $k$. Let*

$a, b \approx 0$. Using (5.6) for $k \in {}^*\mathbb{N} \setminus \mathbb{N}$, we obtain

$$\left(|A_h|h^d\right) = \left(2^k + a\right)\left(\frac{1}{3^k} + b\right)^d$$

$$\leq \left(2^k + a\right)\left(\left(\frac{1}{3^k}\right)^d + b^d\right)$$

$$= 2^k \left(\frac{1}{3^k}\right)^d + a\left(\frac{1}{3^k}\right)^d + 2^k b^d + ab^d$$

$$\approx 2^k \left(\frac{1}{3^k}\right)^d.$$

All the other terms equal $0$ by section 2.6. Hence, the critical point is when $d = \log 2/\log 3$.

Let's try a different variation of the Cantor set.

**Example 5.4.5 (Cantor Set II).** *This is similar to the Cantor set, but instead of dividing the unit interval into 3 equal parts, we will divide it into 5 equal subintervals and delete the second and the fourth. We call this the first iteration and denote it as $C'_1$. So $C'_1$ is the set*

$$\left[0, \frac{1}{5}\right] \cup \left[\frac{2}{5}, \frac{3}{5}\right] \cup \left[\frac{4}{5}, 1\right].$$

*At $C'_2$ we have the intervals*

$$\bigcup_{i=0}^{2}\left[\frac{2i}{5^2}, \frac{2i+1}{5^2}\right] \cup \bigcup_{i=0}^{2}\left[\frac{10+2i}{5^2}, \frac{11+2i}{5^2}\right] \cup \bigcup_{i=0}^{2}\left[\frac{20+2i}{5^2}, \frac{21+2i}{5^2}\right].$$

*We repeat this iteration and $C'$ is generated by taking the intersection of all the $C'_i$'s.*

*The traditional Hausdorff dimension says*

$$H^s_\delta\left(C'\right) = \left(3^k\right)\left(\frac{1}{5^k}\right)^s$$

*gives $s = \log 3/\log 5$. By (5.6), we replace each connected interval with a*

*dot to obtain*

$$
\begin{array}{ll}
k = 0 & \cdot \\
k = 1 & \cdot \quad \cdot \quad \cdot \\
k = 2 & \cdot\,\cdot\,\cdot \quad \cdot\,\cdot\,\cdot \quad \cdot\,\cdot\,\cdot \\
\vdots & \vdots
\end{array}
$$

*We observe that there are $3^k$ dots out of $5^k$ possible dots. The computation is similar as before, and we obtain*

$$\operatorname{st}\left(|A_h|\, h^d\right) = 3^k \left(\frac{1}{5^k}\right)^d,$$

*which implies that $d = \log 3/\log 5$. This corresponds with the Hausdorff dimension which we have found above which is what we wanted.*

*We return to the standard fractals which we have examined earlier.*

**Example 5.4.6 (Sierpinski Gasket).** *As for the Sierpinski gasket (or the Sierpinski triangle), $S_G$ we now replace each connected surface with a point. So for each triangle of side-length $1/(2^k)$ for each iteration, we replace it with a dot. Hence the nonstandard version of our fractal is*

$$
\begin{array}{ccc}
& & \cdot \\
& \cdot & \cdot\,\cdot \\
\cdot & \cdot\,\cdot & \cdot\quad\cdot \\
& & \cdot\,\cdot\quad\cdot\,\cdot \\
k = 0 & k = 1 & k = 2
\end{array}
$$

*Since the number of dots increases by $3^k$ while the maximum number of dots along an edge is $2^k$, we compute*

$$\operatorname{st}\left(|A_h|\, h^d\right) = 3^k \left(\frac{1}{2^k}\right)^d$$

*which gives us $d = \log 3/\log 2$.*

**Example 5.4.7 (Sierpinski Carpet).** As for the Sierpinski carpet $S_C$ we will substitute each square with a point. We observe the initial and the first iteration in the following way:

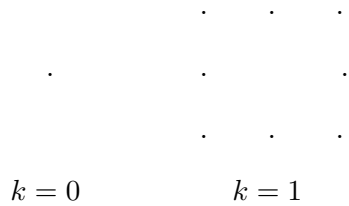

$$k = 0 \qquad\qquad k = 1$$

We have 8 dots when $k = 1$ and we have a maximum of 3 colinear dots along an edge. By our nonstandard equation (5.6),

$$\text{st}\left(|A_h|\, h^d\right) = 8^k \left(\frac{1}{3^k}\right)^d$$

concluding that $d$ is $\log 8 / \log 3$.

**Example 5.4.8 (Quadratic Koch Curve).** We repeat the same process to the quadratic Koch curve $K_Q$ to obtain

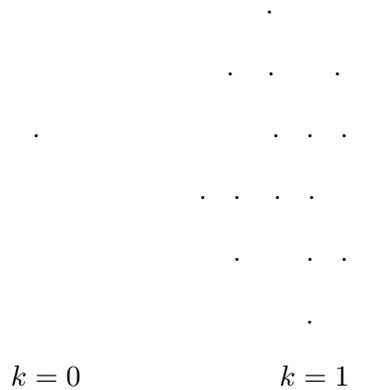

$$k = 0 \qquad\qquad k = 1$$

Since we have $16^k$ dots with $4^k$ possible colinear dots for each iteration $k$,

$$\text{st}\left(|A_h|\, h^d\right) = 16^k \left(\frac{1}{4^k}\right)^d.$$

So $d = 2$.

**Example 5.4.9 (Menger Sponge).** *Finally we replace each cube with a dot for the Menger sponge $M_S$. There is no need to plot the dots since this will only be messy, but the reader should be able to grasp how $M_S$ would look like in the nonstandard world. We will remind the reader that the Menger sponge is a 3-dimensional duplicate of the Sierpinski carpet.*

*Since we have $20^k$ points with $3^k$ points along an edge for each k,*

$$\text{st}\left(|A_h| \, h^d\right) = 20^k \left(\frac{1}{3^k}\right)^d.$$

*And thus $d = \log 20 / \log 3$.*

**Example 5.4.10 (Hyperpyramid).** *Although this is a bit difficult to see, we replace each volume of a shape of a hyperpyramid with a dot resulting in $9^k$ dots with $2^k$ colinear dots. This gives us*

$$\text{st}\left(|A_h| \, h^d\right) = 9^k \left(\frac{1}{2^k}\right)^d,$$

*and so $d = \log 9 / \log 2$.*

## 5.5 Benefits of Nonstandard Hausdorff Measure

Some reasons why we have looked into a nonstandard version of the Hausdorff measure is because for one reason, it introduces "counting" methods that calculate the correct Hausdorff measure. Box-counting methods have been proved to over-estimate or under-estimate the correct dimension of a set, particularly if it is irregular [22]. In fact if a set or its complement is not self-similar, then the box-counting theorem fails since there is no dimension for which the limit converges [33].

Secondly the nonstandard version is much easier to compute than the standard Hausdorff version. This is because the discrete measure of $H^s(A)$ has a fixed $\delta$, not a varying one. Also the "sup" over all $\delta$'s is omitted since the process of taking the supremum over all $\delta$'s has already been applied when choosing our $\delta$. So although we omit this additional step, we continue to obtain the accurate Hausdorff dimension, which is defined for every set.

Thirdly recall the relation between $H^s(A)$ and $c_s\overline{\lambda}(A)$ for nice sets $A$. Because we can easily compute one from another, this gives us various approaches to the same result.

# Chapter 6

# Sets Generated By Sequences

## 6.1 Examples

Suppose $a = \{a_n\}_{n \in \mathbb{N}}$ and $b = \{b_n\}_{n \in \mathbb{N}}$ are sequences where $a_n$'s and $b_n$'s are natural numbers. Let $A_{ab}$ be the set of real numbers whose base 2 representation is

$$0.\overbrace{00\ldots0}^{a_0}\overbrace{xx\ldots x}^{b_0}\overbrace{00\ldots0}^{a_1}\overbrace{xx\ldots x}^{b_1}0\ldots$$

where each $x$ is chosen independently as 0 or 1 and there are $a_0$ 0's in the first block, $b_0$ $x$'s in the second, $a_1$ 0's in the third block, $b_1$ $x$'s in the fourth block, and so on. Then $A = A_{ab}$, which is contained in $[0, 1]$, should have a dimension between 0 and 1. We shall consider an example and verify its Hausdorff dimension.

**Example 6.1.1.** *Consider the sequences $a = b = \{1, 2, 4, 8, 16, 32, \ldots\}$ which define the set*

$$A = \{.0x00xx0000xxxx00000000xxxx\ldots : x = 0, 1\}.$$

*To find the Hausdorff dimension of this set we will cover A in such a way to minimize the overlapping of the covers and minimize the difference between the measure of the set and its covering. Let* $\delta_m = 2^{-m}$, *and denote* $I_i^m = [i2^{-m}, (i+1)2^{-m})$ *which has* $\text{diam}(I_i^m) = \delta_m$. *For each m a certain subset of the* $I_i^m$*'s cover A. We can see how many are needed by observing the first m digits after the binary point.*

*Let* $U_j^m$ *range over those* $I_i^m$*'s that meet A, i.e., the collection* $\{U_j^m\}$ *is a covering of A. The union of* $U_j^m$*'s gives*

$$\bigcup_j U_j^m = \bigcup_i I_i^m \cap A.$$

*Fixing* $\delta = \delta_m$ *we observe the following*

|  | the number of sets in the cover |
|---|---|
| $m = 1$ | 1 |
| $m = 2$ | 2 |
| $m = 3$ | 2 |
| $m = 4$ | 2 |
| $m = 5$ | 4 |
| $m = 6$ | 8 |
| $\vdots$ | $\vdots$ |
| $m = k$ | $2^{\#(x's)}$ |

*where we obtain* $\#(x's)$ *by looking at the first m binary points. For each additional x, number of sets increases by a factor of 2 and for each additional 0, number of sets in our covering does not increase. Hence we obtain*

$$H_\delta^s(A) = \inf \sum_i \text{diam}(U_i)^s$$

$$\leq \sum_j \text{diam}(U_j^m)^s$$

$$= \#(\{U_j^m\}_j) \times (2^{-m})^s.$$

We write $\#\left(\{U_j^m\}_j\right)$ to denote the number of $U_j^m$'s that which has diameter $2^{-m}$.

Since we want to compute $H^s(A)$, it suffices to compute $H^s_{\delta_n}(A)$ for some particular diameter $\delta_n$ from a sequence of $\delta_i$'s. It is because $H^s_\delta$ is monotonically nonincreasing as $\delta$ ranges over the reals, and our definition of $H^s(A)$ takes $\delta \to 0$. Hence, we choose $\delta = 2^{-m}$ to make our computation as simple as possible.

Now, by the definition of $H^s(A) = \lim_{\delta \to 0} H^s_\delta(A) = \sup_{\delta > 0} H^s_\delta(A)$, we want the biggest possible value that which can be obtained when maximizing the number of sets in the covering of $A$ and minimizing $\delta$ as much as possible.

So we let the $m^{th}$ digit of our $\delta$ follow after a string of $x$'s. Its immediately preceding digit is an $x$ while its succeeding digit is a $0$. The location for the possible $m$'s is indicated by the bars

$$A = \{.0x|00xx|0000xxxx|00000000xxxxxxxx|00\ldots : x = 0,\ 1\}.$$

So suppose we have $\#(x\text{'s}) = 1 + 2 + \ldots + 2^{k-2} = 2^{k-1} - 1$ and $\delta = 2^{-(2^k-2)}$. Then we compute

$$H^s_\delta(A) \leq 2^{\left(2^{k-1}-1\right)} \times \left(2^{-\left(2^k-2\right)}\right)^s.$$

Hence

$$H^s(A) = \lim_{\delta \to 0} \frac{2^{\left(2^{k-1}-1\right)}}{2^{s\left(2^k-2\right)}} = \lim_{k \to \infty} 2^{\left(2^{k-1}-1\right)-s\left(2^k-2\right)}$$
$$= 2^{2s-1} \lim_{k \to \infty} 2^{2^{k-1}(1-2s)}.$$

This implies if $s < \frac{1}{2}$ then $H^s(A) \to \infty$. If $s > \frac{1}{2}$, then $H^s(A)$ tends to $0$. Hence the critical point is when $s = \frac{1}{2}$ for this specific $\delta = 2^{-m}$ where $m = \left(2^k - 2\right)$.

Next we use the definition of Hausdorff measure to estimate $\dim_H(A)$. From

$$H^s_\delta(A) \leq 2^{2s-1} \times 2^{2^{k-1}(1-2s)},$$

we obtain $H^s_\delta(A) \to 0$ as $\delta \to 0$ for $s > \frac{1}{2}$. So $s$ must be $\leq \frac{1}{2}$.

Now by a different choice of $\delta$ and similar test intervals, the same argument can show $\dim_H(A) \leq \frac{1}{3}$. We choose our diameter to succeed a chain of zeros, where possible cut-off points for $m$ are indicated by the bars

$$A = \{.0|x00|xx0000|xxxx00000000|xxxx\ldots : x = 0,\ 1\}.$$

So again, we'll fix $\delta = 2^{-m}$, but $m$ changes to $2(2^k - 1) + 2^k = 3 \times 2^k - 2$. Then

$$H^s(A) = \lim_{\delta \to 0} \frac{2^{(2^k - 1)}}{2^{s(3 \times 2^k - 2)}}$$

$$= \lim_{k \to \infty} 2^{2s-1} \times 2^{2^k(1-3s)}$$

$$= 2^{2s-1} \lim_{k \to \infty} 2^{2^k(1-3s)}.$$

Hence the critical point is when $s = \frac{1}{3}$. We also see that $\dim_H(A) \leq \frac{1}{3}$ since we have $H^s(A) \to \infty$ as $k \to \infty$ for $s \leq \frac{1}{3}$.

Now we must prove that the inequality goes in the other direction, i.e., $s \geq \frac{1}{3}$. If we can prove this, then the Hausdorff dimension of $A$ is $1/3$. But if this doesn't work, then the Hausdorff dimension is less than $1/3$ and we must use other means to find its dimension.

So it suffices to show that for some collection of intervals covering $A$ and for any $\epsilon > 0$

$$\sum_i \mathrm{diam}(U_i)^s \geq (1 + \epsilon)^s.$$

So suppose $\{U_i\}$ is a collection of intervals covering $A$. First I will show $A$ is closed. Let's rewrite $A$ as the intersection of the following sets

$$\{.0xxxxx\ldots : x = 0,\ 1\}$$
$$\cap \{.0x00xx\ldots : x = 0,\ 1\}$$
$$\cap \{.0x00xx0000xxxx\ldots : x = 0,\ 1\}$$
$$\cap \ \ldots$$

As we can see from above, the first set is a closed interval $\left[0, \frac{1}{2}\right]$, the second set is $\left[0, \frac{1}{16}\right] \cup \left[\frac{1}{4}, \frac{5}{16}\right]$ and so forth. Since arbitary intersection of closed sets is closed $A$ is closed. And since our set is closed and bounded it is compact. So $A$ has a finite subcover. We will denote the finite subcover as $\{V_j\}_{j=1}^m$.

For each $V_j$ replace with mutually disjoint (except at the endpoints) $k_j$ sets in the construction of $A$. All $k_j$ intervals $I_1^j$, $I_2^j, \ldots I_{k_j}^j$ are of the same size that depend on the number $s$ and the construction of $A$.

I will show for each $j$

$$V_j \cap A \subseteq \bigcup_{i=1}^{k_j} I_i^j \cap A,$$

and

$$\sum_{j=1}^m \operatorname{diam}(V_j)^s \geq \sum_{i,\,j} \operatorname{diam}(I_i^j)^s.$$

If $V_j \cap A = \varnothing$, we will let $I_1^j = I_2^j = \ldots = I_{k_j}^j = \varnothing$.

If $V_j \cap A \neq \varnothing$, let $a = \inf(V_j \cap A)$ and $b = \sup(V_j \cap A)$. Since $A$ is closed, the elements $a$ and $b$ are in $A$.

There are few cases to consider.

1. If $a = b$ then we have a singleton point and we choose only one set $I_1$ of diameter $\epsilon > 0$ centered around $a$.

2. If $a \neq b$ and $I_i^j \cap I_k^j = \{c\}$ for some $i \neq k$ then standard intervals $I_i^j$ and $I_k^j$ intersect only in their endpoint. Without loss of generality, assume $\inf I_i^j \leq \inf I_k^j$ and $\sup I_i^j \leq \sup I_k^j$. Then since they intersect at exactly one point $\sup I_i^j = \inf I_k^j = c$ and

$$\operatorname{diam}(I_i^j \cup I_k^j) = 2\operatorname{diam}(I_i^j).$$

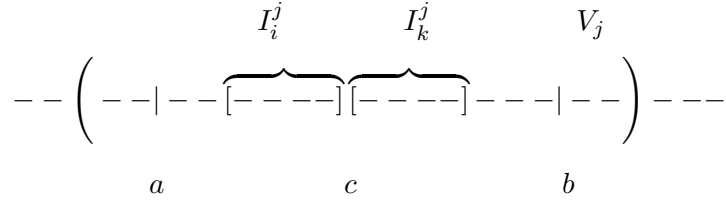

3. For $a \neq b$ and $c_1, c_2 \in I_i^j \cap I_k^j$ where $i \neq k$ and $a \lneq c_1 \lneq c_2 \lneq b$

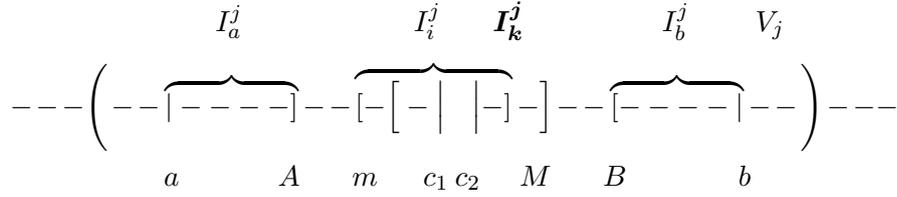

let $m = \inf(I_i^j, I_k^j)$ and $M = \sup(I_i^j, I_k^j)$. Since $a, b \in A$ there are some intervals $I_a^j$ and $I_b^j$ so that $a \in I_a^j$ and $b \in I_b^j$.

If $a = \inf I_b^j$ or $b = \sup I_a^j$ then we have exactly one interval $I_a^j = I_b^j$ of diameter $b - a$. So suppose not.

Let $A = \sup I_a^j$ and $B = \inf I_b^j$. By our choice of standard intervals $I_1^j, I_2^j, \ldots, I_{k_j}^j$, they have equal diameter. So here they have diameter of $A - a = b - B = M - m$. By our construction of $A$ and since $I_i^j$ and $I_k^j$ intersect in at least two disjoint points, $I_i^j = I_k^j$.

Hence for our chosen $I_1^j, \ldots, I_{k_j}^j$ we will assume each $I_i^j \cap I_k^j = \emptyset$ for every $i \neq k$, or they intersect at most one or two points, which is at the endpoints.

Now let
$$K = V_j \setminus \bigcup_{i=1}^{k_j} I_i^j.$$
Since the unit interval is halving at each iteration
$$\operatorname{diam}\left(I_1^j\right), \operatorname{diam}\left(I_2^j\right), \ldots, \operatorname{diam}\left(I_{k_j}^j\right) \leq \operatorname{diam}(K).$$
Moreover we chose our $I_i^j$'s so that $K \cap I_i^j = \emptyset$ for each $i$. Hence by our construction of $I_i^j$'s
$$V_j \cap A \subseteq \bigcup_{i=1}^{k_j} I_i^j \cap A$$

*for each $j$, and since each diameter of $I_i^j$'s is less than or equal to the diameter of $K$,*

$$\mathrm{diam}(V_j)^s = \mathrm{diam}\left(\bigcup_i \left(I_i^j \cup K\right)\right)^s \geq \sum_{i=1}^{k_j} \mathrm{diam}(I_i^j)^s,$$

*or*

$$\sum_j \mathrm{diam}(V_j)^s \geq \sum_{i,\,j} \mathrm{diam}\left(I_i^j\right)^s.$$

*This shows that the right-hand side doesn't increase the sum. Finally we need to show $\sum_{i,\,j} \mathrm{diam}\left(I_i^j\right)^s = 1$ for $s \geq \frac{1}{3}$.*

*Since the $I_i^j$'s are in the construction of $A$, they are of "standard" size. We have shown earlier that we need $2^{(2^{k-1}-1)}$ sets if the diameter is of size $2^{-(2(2^{k-1}-1))} = 2^{-(2^k-2)}$. And if our diameter is of length $2^{-(3\times 2^{k-1}-2)}$, then we still need $2^{(2^{k-1}-1)}$ sets to cover $A$. Either way, $\sum_i \mathrm{diam}(U_i)^s = 1$, so for $s \geq \frac{1}{3}$, $\sum_{i,\,j} \mathrm{diam}(I_i^j)^s = 1$. So we conclude*

$$\sum_{j=1}^m \mathrm{diam}(V_j)^s \geq 1.$$

*This shows that $s \geq \frac{1}{3}$. Hence, the Hausdorff dimension $s$ for $A$ is $\frac{1}{3}$.*

Now let's change the sequences and examine the relationship between our generated set and its Hausdorff dimension. From this point, we will always choose our $\delta$-covering to succeed a sequence of zeros since we have proved that this gives us the correct Hausdorff dimension.

**Example 6.1.2 (Geometric Sequences).** *Let $n \in \mathbb{N}$ be strictly positive and let $a$ and $b$ be the geometric sequence $\{n^i\}_{i \in \mathbb{N}}$ generating the set*

$$A = \{.\underbrace{0\ldots 0}_{n^0}\underbrace{x\ldots x}_{n^0}\underbrace{0\ldots 0}_{n^1}\underbrace{x\ldots x}_{n^1}\ldots\}.$$

*Let's use $|x|$ to represent the number of distinct elements in $x$. So for $n = 1$, we need $|x|^{(k-1)}$ sets of diameter $|x|^{-(2k-1)}$ to cover $A$. For $n > 1$, we need*

$$|x|^{\left(\frac{1-n^{k-1}}{1-n}\right)}$$

sets of diameter
$$|x|^{-\left(2\left(\frac{1-n^{k-1}}{1-n}\right)+n^{k-1}\right)}$$
to cover A. Hence we obtain a table summarizing the Hausdorff dimension of $A$ as $n$ and $|x|$ vary over the naturals.

Table of $\dim_H(A)$ for Geometric Sequences

| $n^i$ | $1^i$ | $2^i$ | $3^i$ | $4^i$ | $5^i$ | ... | $16^i$ | ... |
|---|---|---|---|---|---|---|---|---|
| $x = 0, 1$ | $\frac{1}{2}$ | $\frac{1}{3}$ | $\frac{1}{4}$ | $\frac{1}{5}$ | $\frac{1}{6}$ | ... | $\frac{1}{17}$ | ... |
| $x = 0, 1, 2$ | $\frac{1}{2}$ | $\frac{1}{3}$ | $\frac{1}{4}$ | $\frac{1}{5}$ | $\frac{1}{6}$ | ... | $\frac{1}{17}$ | ... |
| $x = 0, 1, 2, 3$ | $\frac{1}{2}$ | $\frac{1}{3}$ | $\frac{1}{4}$ | $\frac{1}{5}$ | $\frac{1}{6}$ | ... | $\frac{1}{17}$ | ... |
| $x = 0, 1, 2, 3, 4$ | $\frac{1}{2}$ | $\frac{1}{3}$ | $\frac{1}{4}$ | $\frac{1}{5}$ | $\frac{1}{6}$ | ... | $\frac{1}{17}$ | ... |
| $\vdots$ | | | $\vdots$ | | | | | |

So as $n$ increases, the Hausdorff dimension decreases in the rate of $(1+n)^{-1}$ which is independent of our choice for $x$.

Next we summarize the Hausdorff s-dimensional outer measure of $A$.

Table of $H^s(A)$ for Geometric Sequences

| $n^i$ | $1^i$ | $2^i$ | $3^i$ | $4^i$ | $5^i$ | ... | $16^i$ | ... |
|---|---|---|---|---|---|---|---|---|
| $x = 0, 1$ | $\frac{1}{\sqrt[2]{2}}$ | $\frac{1}{\sqrt[3]{2}}$ | $\frac{1}{\sqrt[4]{2}}$ | $\frac{1}{\sqrt[5]{2}}$ | $\frac{1}{\sqrt[6]{2}}$ | ... | $\frac{1}{\sqrt[17]{2}}$ | ... |
| $x = 0, 1, 2$ | $\frac{1}{\sqrt[2]{3}}$ | $\frac{1}{\sqrt[3]{3}}$ | $\frac{1}{\sqrt[4]{3}}$ | $\frac{1}{\sqrt[5]{3}}$ | $\frac{1}{\sqrt[6]{3}}$ | ... | $\frac{1}{\sqrt[17]{3}}$ | ... |
| $x = 0, 1, 2, 3$ | $\frac{1}{\sqrt[2]{4}}$ | $\frac{1}{\sqrt[3]{4}}$ | $\frac{1}{\sqrt[4]{4}}$ | $\frac{1}{\sqrt[5]{4}}$ | $\frac{1}{\sqrt[6]{4}}$ | ... | $\frac{1}{\sqrt[17]{4}}$ | ... |
| $x = 0, 1, 2, 3, 4$ | $\frac{1}{\sqrt[2]{5}}$ | $\frac{1}{\sqrt[3]{5}}$ | $\frac{1}{\sqrt[4]{5}}$ | $\frac{1}{\sqrt[5]{5}}$ | $\frac{1}{\sqrt[6]{5}}$ | ... | $\frac{1}{\sqrt[17]{5}}$ | ... |
| $\vdots$ | | | $\vdots$ | | | | | |

So the $H^s(A)$ increases as $n$ increases which corresponds with the Hausdorff measure definition. It says $H^s_\delta(A)$ increases as $\delta$ decreases.

**Example 6.1.3.** *Next suppose a and b are monotonically increasing arithmetic sequences of the form*

$$1 + 0d,\ 1 + 1d,\ 1 + 2d,\ 1 + 3d, \ldots$$

*where d is a nonnegative integer. Then as d and $|x|$ vary over the naturals, the Hausdorff dimension of the set is summarized as before.*

Table of $\dim_H(A)$ for Arithmetic Sequences

| $d$ | 0 | 1 | 2 | 3 | 4 | ... | $d'$ | ... |
|---|---|---|---|---|---|---|---|---|
| $x = 0, 1$ | $\frac{1}{2}$ | $\frac{1}{2}$ | $\frac{1}{2}$ | $\frac{1}{2}$ | $\frac{1}{2}$ | ... | $\frac{1}{2}$ | ... |
| $x = 0, 1, 2$ | $\frac{1}{2}$ | $\frac{1}{2}$ | $\frac{1}{2}$ | $\frac{1}{2}$ | $\frac{1}{2}$ | ... | $\frac{1}{2}$ | ... |
| $x = 0, 1, 2, 3$ | $\frac{1}{2}$ | $\frac{1}{2}$ | $\frac{1}{2}$ | $\frac{1}{2}$ | $\frac{1}{2}$ | ... | $\frac{1}{2}$ | ... |
| $x = 0, 1, 2, 3, 4$ | $\frac{1}{2}$ | $\frac{1}{2}$ | $\frac{1}{2}$ | $\frac{1}{2}$ | $\frac{1}{2}$ | ... | $\frac{1}{2}$ | ... |
| $\vdots$ | | | | | | | $\vdots$ | |

*Since the Hausdorff dimension of A is not affected by the arithmetic growth rate, this result tells us that arithmetic sequences do not increase fast enough. In fact, sets generated by arithmetic sequences must be in a class of their own resulting in identical patterns. This table summarizes the Hausdorff s-dimensional outer measure.*

Table of $H^s(A)$ for Arithmetic Sequences

| $d$ | 0 | 1 | 2 | 3 | 4 | ... | $d'$ | ... |
|---|---|---|---|---|---|---|---|---|
| $x = 0, 1$ | $\frac{1}{\sqrt{2}}$ | $\frac{1}{\sqrt{2}}$ | $\frac{1}{\sqrt{2}}$ | $\frac{1}{\sqrt{2}}$ | $\frac{1}{\sqrt{2}}$ | ... | $\frac{1}{\sqrt{2}}$ | ... |
| $x = 0, 1, 2$ | $\frac{1}{\sqrt{3}}$ | $\frac{1}{\sqrt{3}}$ | $\frac{1}{\sqrt{3}}$ | $\frac{1}{\sqrt{3}}$ | $\frac{1}{\sqrt{3}}$ | ... | $\frac{1}{\sqrt{3}}$ | ... |
| $x = 0, 1, 2, 3$ | $\frac{1}{\sqrt{4}}$ | $\frac{1}{\sqrt{4}}$ | $\frac{1}{\sqrt{4}}$ | $\frac{1}{\sqrt{4}}$ | $\frac{1}{\sqrt{4}}$ | ... | $\frac{1}{\sqrt{4}}$ | ... |
| $x = 0, 1, 2, 3, 4$ | $\frac{1}{\sqrt{5}}$ | $\frac{1}{\sqrt{5}}$ | $\frac{1}{\sqrt{5}}$ | $\frac{1}{\sqrt{5}}$ | $\frac{1}{\sqrt{5}}$ | ... | $\frac{1}{\sqrt{5}}$ | ... |
| $\vdots$ | | | | | | | $\vdots$ | |

*Hence variations of arithmetic sequences do not increase fast enough to affect $H^s(A)$. It only depends on the size of $x$.*

Next, we will consider a power tower sequence

$$a = b = \{2^0,\ 2^1,\ 2^2,\ 2^{2^2},\ 2^{2^{2^2}}, \ldots\}.$$

In order to analyze the set generated by these sequences, we will need to first prove a lemma.

**Lemma 6.1.4.** *For all $x \geq 7$, $2^x > 2^{x-1} + x^2$.*

*Proof.* Let's prove by induction. Let $x = 7$. Then $2^7 > 2^{7-1} + 7^2$. So assume our induction hypothesis that $2^x > 2^{x-1} + x^2$ is true for some $x$. Then

$$2^{x+1} = 2 \times 2^x > 2\left(2^{x-1} + x^2\right)$$
$$= 2^x + 2x^2$$
$$> 2^x + (x+1)^2.$$

Hence $2^x > 2^{x-1} + x^2$ for each $x \geq 7$. $\square$

This lemma says $2^x$ increases faster than $\left(2^{x-1} + x^2\right)$ for any $x \geq 7$. Now we introduce a definition.

**Definition 6.1.5.** *Define $2_k = 2^{2^{2^{\cdots 2}}}$ where 2 is raised to the power of two $(k-1)$-times. Our notation $2_k$ is called the power tower of 2.*

**Proposition 6.1.6.** *Let $x_i = (x_{i-1})^2$ for each $i \in \mathbb{N}$, and we let $x_0 = 2^0$ and $x_1 = 2^1$. Then for any $\epsilon > 0$ there is some $N_0$ so that for all $k \geq N_0$*

$$\sum_{i=0}^{k} x_i \leq (1+\epsilon)\, x_k.$$

So, what we want to show is that $1 + 2 + \ldots + 2_{k-1} + 2_k \approx 2_k$ for all large $k$.

*Proof.* Let's put $1 + 2 + \ldots + 2_{k-1}$ to the other side to obtain

$$2_k \approx 2_k - (1 + 2 + \ldots + 2_{k-1}).$$

Now consider $2_k - 2_{k-1} - \ldots - 1$.

$$\begin{aligned}
2_k - 2_{k-1} - \ldots - 1 &> 2_k - (2_{k-1})^2 \\
&= 2^{(2_{k-1})} - (2_{k-1})^2 \\
&> 2^{(2_{k-1})-1} \\
&\gg 2_{k-1} = \frac{1}{2} 2_k
\end{aligned}$$

where we have used Lemma 6.1.4 to obtain the inequality on the third line. Hence $\sum_{i=0}^{k} 2_i \approx 2_k$. $\square$

Now we analyze the power tower of 2.

**Example 6.1.7.** *Let a and b be the sequence*

$$\{2^0,\ 2^1,\ 2^2,\ 2^{2^2},\ 2^{2^{2^2}},\ldots\}$$

*which generates*

$$P = \{.0x00xx0000xxxx\ldots : x = 0,\ 1\}.$$

*Since the diameter of the covering of P follows a chain of zeros, we let our diameter be at the $\left(2\left(1 + 2 + 2^2 + \ldots + 2_k\right) + 2_{k+1}\right)^{th}$-digit which is indicated by the bar*

$$P = \{.0x00xx\ldots \underbrace{xx\ldots x}_{2_k} \underbrace{00\ldots 0}_{2_{k+1}} |xx\ldots : x = 0,\ 1\}.$$

*Use Proposition 6.1.6 to obtain $2\left(1 + 2 + 2^2 + \ldots + 2_k\right) + 2_{k+1} \approx 2_{k+1}$, or*

$$2^{-\left[2\left(1+2+2^2+\ldots+2_k\right)+2_{k+1}\right]} \approx 2^{-[2_{k+1}]}.$$

*We used the same proposition once more to conclude that we need $2^{[2_k]}$-many sets with this diameter to cover $P$. So*

$$H_\delta^s(P) \leq 2^{2_k} \times 2^{-2_{k+1} \times s}$$
$$= 2^{2_k} \times 2^{-2^{2_k} \times s}$$
$$= 2^{2^{(2_k-1)} - 2^{(2_k)} \times 2^{(\log_2 s)}}$$
$$= 2^{2^{(2_k-1)} - 2^{(2_k + \log_2 s)}}.$$

*So if we assume $s > 0$, we have $2^{-2^{(2_k + \log_2 s)}}$ growing much faster than $2^{2^{(2_k-1)}}$. This means the right-hand side tends to $-\infty$. Since there isn't any positive $s$ for which the limit converges, $s = 0$ [33]. Hence we conclude $H^s(P) \to \infty$ as $k \to \infty$ for $s = 0$.*

*This is an example of a set with Hausdorff dimension $0$ and Hausdorff $s$-dimensional outer measure $\infty$.*

The power tower of 2 led us to believe that there might be at least one another sequence between $\{2^n\}$ and $\{2_n\}$ so that this sequence generates a set whose critical Hausdorff dimension is 0. Hence this leads to our final analysis of these particular type of sets.

Before we go on about to find these specific sequences, we list some conditions on our sequences. First where $a = \{a_n\}_{n \in \mathbb{N}}$ and $b = \{b_n\}_{n \in \mathbb{N}}$, we assume $a_n = b_n$ for almost all $n$. Secondly, we want $a$ and $b$ to be nondecreasing sequences. And lastly, if they increase fast enough, then the zeros must dominate since the correct $\dim_H(A)$ is the lower box-counting dimension as we have shown earlier.

What we mean by the last requirement is the following: let $A$ be the set

$$0.\overbrace{0\ldots0}^{a_0}\overbrace{x\ldots x}^{a_0}\overbrace{0\ldots0}^{a_2}\ldots\ldots\overbrace{x\ldots x}^{a_{n-1}}\overbrace{0\ldots0}^{a_n}|\ldots$$

where $x$ represents $0, 1$ as before. We obtain our $\delta$-cover where its cut-off point is after some $a_n$-zeros. Then $\#(x's) = \sum_{i=0}^{n-1} a_i$ and $\#(0's) = \sum_{i=0}^{n-1}$

$a_i + a_n$. Hence our diameter is of size

$$2^{-(\#(x's)+\#(0's))} = 2^{-\left(2\sum_{i=0}^{n-1} a_i + a_n\right)}.$$

We observe that we want $a_n$ to be as big as it can get so that it dominates all the $a_i$'s. In other words,

**Conjecture 6.1.8.** *Let $a = \{a_n\}_{n\in\mathbb{N}}$ be a sequence of natural numbers. Then if $a_n > F\left(\sum_{i=0}^{n-1} a_i\right)$ for infinitely-many $n$ and a constant $F$, then $\dim_H(A) = 0$.*

So computing the Hausdorff measure,

$$H_\delta^s(A) \leq \left(2^{\sum_{i=0}^{n-1} a_i}\right)\left(2^{-\left(2\sum_{i=0}^{n-1} a_i + a_n\right)}\right)^s$$

$$\leq 2^{(1-2s)\sum_{i=0}^{n-1} a_i - sa_n}$$

implies we obtain the critical value for $s$ when

$$(1-2s)\sum_{i=0}^{n-1} a_i - sa_n \tag{6.1}$$

is a positive finite real number.

First, we know $s$ must be less than $1/2$ for if we suppose $s \geq 1/2$, or namely, $s = 1/2$, then we simplify (6.1) to $-(1/2)a_n$ which diverges to $-\infty$ as $n \to \infty$ since $\{a_n\}$ is unbounded. So $0 \leq H_\delta^{1/2}(A) \leq 2^{-(1/2)a_n} \to 0$ implies $H_\delta^{1/2}(A) = 0$. This means $s = 1/2$ is too big. So $\dim_H(A)$ must be $< 1/2$. Also by the definition $\dim_H(A) = \inf\{s\colon H^s(A) = 0\}$, we want the smallest $s$ so that this holds true.

Secondly, we know $s \geq 0$ but can we build another requirement so that we may conclude $s = 0$?

**Conjecture 6.1.9.** *Continuing from above, if $\{a_n\}_{n\in\mathbb{N}}$ is nondecreasing sequence and for any $s > 0$,*

$$a_n - \left(\frac{1-2s}{s}\sum_{i=0}^{n-1} a_i\right) \to +\infty \tag{6.2}$$

*if, and only if, $\dim_H(A) = 0$.*

In simple terms, for any $K > 0$ we can rewrite (6.2) as

$$a_n - K \sum_{i=0}^{n-1} a_i \to +\infty. \tag{6.3}$$

Hence it is sufficient to find a sequence so that (6.2) holds.

**Example 6.1.10.** *Consider the sequence $a = \{2^{2^n}\}_{n \in \mathbb{N}}$. We know*

$$\sum_{i=0}^{n-1} a_i = 2^{2^0} + 2^{2^1} + \ldots + 2^{2^{n-1}} \leq 2^{2^{n-1}+1}.$$

*So for some fixed $K > 0$, we use the sequence $\{2^{2^n}\}_{n \in \mathbb{N}}$ to compute (6.3)*

$$a_n - K \sum_{i=0}^{n-1} a_i = 2^{2^n} - K \sum_{i=0}^{n-1} 2^{2^i}$$

$$\geq 2^{2^n} - K 2^{2^{n-1}+1} = 2^{2^{n-1}+2^{n-1}} - K 2^{2^{n-1}+1}$$

$$= 2^{2^{n-1}+(2^{n-1})} - 2^{2^{n-1}+(\log_2 K + 1)}.$$

*Since $2^{n-1}$ grows faster than $\log_2 K + 1$, everything on the last line tends to infinity, which also pushes the top left-hand side to infinity satisfying (6.3). And so this sequence generates a set with Hausdorff dimension 0.*

**Example 6.1.11.** *Let's analyze another sequence $\{3^{3^n}\}_{n \in \mathbb{N}}$. Then the sum of its first $(n-1)$ terms is*

$$\sum_{i=0}^{n-1} a_i = 3^{3^0} + 3^{3^1} + \ldots + 3^{3^{n-1}} \leq \frac{1 - 3^{3^{n-1}+1}}{1 - 3}.$$

*So computing the left-hand side of (6.3),*

$$a_n - K \sum_{i=0}^{n-1} a_i = 3^{3^n} - K \sum_{i=0}^{n-1} 3^{3^i}$$

$$\geq 3^{3^n} - K \left( \frac{3^{3^{n-1}+1} - 1}{2} \right)$$

$$\geq 3^{3^n} - \frac{K}{2} \left( 3^{3^{n-1}+1} \right)$$

$$= 3^{3^n} - 3^{3^{n-1}+1+\log_3(K/2)}$$

$$= 3^{3^{n-1}+(3^{n-1}+3^{n-1})} - 3^{3^{n-1}+(\log_3(K/2)+1)}.$$

Similar as before, $(3^{n-1} + 3^{n-1})$ is grows faster than $(\log_3(K/2) + 1)$ which means property (6.3) is satisfied, and so $\dim_H(A) = 0$.

Now what about nondecreasing sequences that satisfy

$$a_n = \left(\sum_{i=0}^{n-1} a_i\right)^2 ? \tag{6.4}$$

Would a sequence that which meets this condition also fulfill (6.3)? If we let $a_0 = 1$, then we can easily compute the rest of the $a_i$'s: $a_1 = 1$, $a_2 = (2)^2 = 4$, $a_3 = (6)^2 = 36$, $a_4 = (42)^2 = 1764, \ldots$. But we want to approximate its growth rate, not simply list out as many $a_i$'s as much as possible, so that we may determine other sequences that fulfill (6.3).

So first of all, we know $\{3^{3^n}\}_{n \in \mathbb{N}}$ meets this condition since

$$\left(\sum_{i=0}^{n-1} a_i\right)^2 = \left(3^{3^0} + 3^{3^1} + \ldots + 3^{3^{n-1}}\right)^2$$

$$\leq \left(\frac{1}{2}\left(3^{3^{n-1}+1} - 1\right)\right)^2$$

$$< \frac{1}{4}\left(3^{2 \times 3^{n-1}+2}\right) \leq 3^{3^n} = a_n.$$

But as for $\{2^{2^n}\}_{n \in \mathbb{N}}$, the condition (6.4) is not fulfilled as

$$\left(\sum_{i=0}^{n-1} a_i\right)^2 \leq \left(2^{2^{n-1}+1}\right)^2$$

$$= 2^{2^n+2} \not\leq a_n = 2^{2^n} = a_n.$$

This implies that $\{2^{2^n}\}_{n \in \mathbb{N}}$ doesn't grow in the rate of (6.4). But nevertheless, we may still use this property to find other conditions and sequences that satisfy (6.3). This leads us to the next remark.

**Remark 6.1.12.** *We have a large list of conditions we may use to generate sequences $\{a_n\}_{n \in \mathbb{N}}$ that fulfill (6.3). Such includes:*

1. $a_n = \left(\sum_{i=0}^{n-1} a_i\right) \times \log\left(\sum_{i=0}^{n-1} a_i\right)$;

2. $a_n = \left(\sum_{i=0}^{n-1} a_i\right) \times \log\log\left(\sum_{i=0}^{n-1} a_i\right)$;

3. $a_n = \left(\sum_{i=0}^{n-1} a_i\right) \times f\left(\sum_{i=0}^{n-1} a_i\right)$ where $f$ is an increasing function.

We may also consider sequences in between $2^{2^n}$ and $3^{3^n}$ such as:

1. $a_n = 2^{2^{(\log_2 n)^2}}$;

2. $a_n = 2^{2^{(\log_2 n)(\log_2 \log_2 n)}}$.

**Remark 6.1.13.** If $a = \{b^n\}_{n \in \mathbb{N}}$ where $b \geq 2$, then

$$\sum_{i=0}^{n-1} a_i = \frac{b^n - 1}{b-1} = \frac{1}{b-1} a_n - \frac{1}{b-1}$$
$$= \frac{1}{b-1} a_n - \epsilon.$$

This shows $a_n$ is not much bigger than the sum of the $a_i$'s, so this sequence does not meet the condition (6.3).

## 6.2 Relationship with the Cantor Set

Let's recall how one can generate the Cantor set $C$. Let $C_0$ be the unit interval $[0, 1]$. Break the middle $1/3$ of the interval and let $C_1$ be the intervals $[0, 1/3] \cup [2/3, 1]$. Again take away the middle $1/3$ of each interval and let $C_2$ be $[0, 1/9] \cup [2/9, 3/9] \cup [6/9, 7/9] \cup [8/9, 1]$. Repeat this cutting of its previous intervals and take the intersection of all $C_i$'s to obtain $C$.

Now we might ask ourselves how our particular set is related to the Cantor set? Both our set and the Cantor set are totally disconnected [8], that is, for any distinct $x$ and $y$ in the set, there is some $z$ between them that is not in the set. Also both of these sets are first category, which means the set is a union of countably-many nowhere dense subsets, and the complement of our set and the complement of the Cantor set are second category since they are not first category.

Moreover, both of the complements are open and dense (and so our sets are closed). They are regarded as very large in Baire category theory. Baire category theorem says any complete metric space is a Baire space, and since $[0, 1] \setminus C$ is in a complete metric space, and so are the sets generated by sequences, it is also in a Baire space.

Since the Cantor set has Lebesgue measure $\leq 2^i/3^i$ at each $C_i$, this implies
$$\overline{\lambda}(C) = \overline{\lambda}\left(\bigcap_{i \in \mathbb{N}} C_i\right) = \lim_{i \to \infty} \overline{\lambda}(C_i) \leq \lim_{i \to \infty} \left(\frac{2}{3}\right)^i = 0.$$
So since the outer measure of $C$ is 0, the Cantor set is a "big" set with Lebegue outer measure 0. But when using the definition of Hausdorff outer measure, the Cantor set has measure 1 with $\dim_H(C) = \log 2/\log 3$. Similarly, our particular set has Lebesgue measure 0, but a positive Hausdorff outer measure with $0 \leq \dim_H(A) < 1$. We note that the Hausdorff dimension equals 0 if the sequences used to generate our set satisfy (6.3).

Also local density for both the Cantor set and our set fails to exist at almost every point since they are totally disconnected. So as we can see, the Cantor set has many likeness with our set.

Here, we do a little variation of the Cantor set so that its measure is not zero. Let $a_0 = 1$ and let $a_0, a_1, \ldots$ be a strictly decreasing sequence of points converging to some positive number $a \in (0, 1)$. Let $a_i$ be the measure of the set at its $i^{th}$-iteration. So we have defined a new "Cantor set" $C^a$ so that each $C_i^a$ is composed of the following sets

$$C_0^a = [0, a_0] = [0, 1]$$
$$C_1^a = \left[0, \frac{a_1}{2}\right] \cup \left[1 - \frac{a_1}{2}, 1\right]$$
$$C_2^a = \left[0, \frac{a_2}{4}\right] \cup \left[\frac{a_1}{2} - \frac{a_2}{4}, \frac{a_1}{2}\right] \cup \left[\left(1 - \frac{a_1}{2}\right), \left(1 - \frac{a_1}{2}\right) + \frac{a_2}{4}\right] \cup \left[\left(1 - \frac{a_2}{4}\right), 1\right]$$
$$\vdots \qquad \vdots$$

We take the intersection of all $C_i^a$'s to get $C^a$ and we can see that the outer measure of $C^a$ is $\lim_{i \to \infty} \overline{\lambda}(C_i^a) = a$. So $C^a$ has a positive Lebesgue measure. Similarly with few variations on our particular set, we are able to generate sets with positive Lebesgue measure.

Another quality of the Cantor set is that it corresponds to paths through the full binary tree. For example, let the two nodes on the second row of the binary tree represent the endpoints 0, 1, which are the only endpoints in $C_0$, the Cantor set's initial iteration. We also let each set of nodes along the same level represent the endpoints of the intervals of the Cantor set at each iteration.

We observe that we can follow any of the paths down the binary tree to obtain any of the points in the Cantor set. Since the infinitely growing binary tree and the Cantor set are of the same size, there is a homeomorphism between the Cantor set and the binary tree. That is, $C \cong \{0, 1\}^{\mathbb{N}}$ when the right-hand side is given the product topology. And the size of the Cantor set is equal to $2^{\text{card}(\mathbb{N})} = 2^{\aleph_0}$.

## 6.3 Box-Counting Dimension

We now move-on to analyze the relationship between our sets generated by sequences and the box-counting dimension. We look at an example by returning to the binary set generated by $a$ and $b$, which are sequences with terms

$$1,\ 2,\ 2^2,\ 2^3,\ 2^4, \ldots.$$

We will use the standard $d_{bc}$ definition to compute the box-counting dimension of the set generated by the above sequence. Recalling some old definitions, the lower and the upper box-counting dimensions are

$$\underline{d}_{bc}(A) = \liminf_{\epsilon \to 0} \frac{-\log N_\epsilon(A)}{\log \epsilon}$$

and

$$\overline{d}_{bc}(A) = \limsup_{\epsilon \to 0} \frac{-\log N_\epsilon(A)}{\log \epsilon}$$

respectively.

So let $\epsilon > 0$ as this will be our step-size. Our $N_\epsilon(A)$ is the number of "boxes" intersecting $A$. There are two cases which must be considered. One of them is that the diameter of our covering comes after the $x$'s and the second case is when the diameter comes after the 0's. Setting our diameter to one of these sizes is sufficient since we eventually take $\epsilon$ to 0. We will not worry about the rate in which $\epsilon$ converges to 0 since we want the infimum of all coverings of $A$. This allows us to choose any $\epsilon$ to find the dimension.

Similar as before, we first choose $\epsilon$ to come after the $x$'s. So $\epsilon = 2^{-m}$ where $m = 2 \times 2^{k-1} - 2$. Then we compute $N_\epsilon(A)$ which is $2^{2^{k-1}-1}$. Thus $\overline{d}_{bc}(A)$ is

$$\begin{aligned}
\limsup_{\epsilon \to \infty} \frac{-\log N_\epsilon(A)}{\log \epsilon} &= \lim_{k \to \infty} \frac{-\log 2^{2^{k-1}-1}}{\log\left(2^{-(2 \times 2^{k-1} - 2)}\right)} \\
&= \lim_{k \to \infty} \frac{(2^{k-1} - 1)\log 2}{(2 \times 2^{k-1} - 2)\log 2} \\
&= \lim_{k \to \infty} \frac{1 - 1/(2^{k-1})}{2 - 2/(2^{k-1})} \approx \frac{1}{2}.
\end{aligned}$$

Next we compute the step-size $\epsilon$ as it decreases in such a way so that it follows each sequence of zeros. So $\epsilon = 2^{-n}$ where $n = 3 \times 2^{k-1} - 2$.

So this gives $N_\epsilon(A)$ as $2^{2^{k-1}-1}$. Hence we obtain

$$\underline{d}_{bc}(A) = \liminf_{\epsilon \to 0} \frac{-\log N_\epsilon(A)}{\log \epsilon} = \lim_{k \to \infty} \frac{-\log 2^{2^{k-1}-1}}{\log\left(2^{-(3 \times 2^{k-1}-2)}\right)}$$
$$= \lim_{k \to \infty} \frac{\left(2^{k-1} - 1\right) \log 2}{(3 \times 2^{k-1} - 2) \log 2}$$
$$= \lim_{k \to \infty} \frac{1 - 1/\left(2^{k-1}\right)}{3 - 2/\left(2^{k-1}\right)} \approx \frac{1}{3}.$$

Hence the inner box-counting dimension equals our Hausdorff dimension of $A$.

This shows that our sets are slightly more irregular than typical fractal examples. If the step-size decreased after each chain of $x$'s, then the covering of our set was too big. But when it was after the zeros, we find the correct dimension of our set.

## 6.4 Self-Similarity

Finally, in this section, we briefly discuss self-similarity, and conclude if our set is self-similar.

According to one paper [33], a set $A$ is self-similar if any subset of it can be scaled to look like the original set by multiplying by an appropriate power of 10. In this paper I will call a set self-similar if "it is built up of scaled down pieces that are geometrically similar to the entire set" [7]. For example, G. Edgar uses the example of the Cantor set to show self-similarity [7]. The Cantor set $C$ can be rewritten as $C_1 \cup C_2$ where $C_1 = C \cap [0, \frac{1}{3}]$ and $C_2 = C \cap [\frac{2}{3}, 1]$ with each $C_1$ and $C_2$ being a scaled down copies of $C$ by a factor of 1/3.

Next definition has been extracted from J. E. Hutchinson's paper [18].

**Definition 6.4.1.** *Let $X$ be a complete metric space, and we will let $f_1$, $f_2$, ..., $f_n : X \to X$ be contractive mappings. Then there is a unique nonempty compact subset $K$ of $X$ so that*

$$K = \bigcup_{i=1}^{n} f_i(K).$$

*The set $K$ is called the invariant, or self-similar, set.*

M. Barnsley first introduced the terminology "iterated function system" or IFS [1] and called the list $(X, f_1, \ldots, f_n)$ by this name. We recall the definition of a contraction.

**Definition 6.4.2.** *Let $f : \mathbb{R}^n \to \mathbb{R}^n$ be a function. Then for any $x$ and $y$ in $\mathbb{R}^n$ and for some real $0 \leq C < 1$, if*

$$|f(x) - f(y)| \leq C |x - y|,$$

*then we call $f$ a contraction map.*

So each contraction $f_i$ has some positive constant $C_i < 1$ so that for any distinct $x$ and $y$, $|f_i(x) - f_i(y)| = C_i |x - y|$. And since the union of

$f_i(K)$'s is similar to $K$, and $f_i(K) \cap f_j(K)$ is almost empty for any $i \neq j$, the Hausdorff dimension $s$ also satisfies $\sum_{i=1}^{n} C_i^s = 1$.

**Example 6.4.3.** *Consider the Sierpinski carpet $S_C$ after two iterations. This means that we have cut-out a square of side-length $1/3$ and $8$ squares of side-length $1/(3^2)$. We zoom-in on the top left corner, which looks like a copy of the Sierpinski carpet after its first iteration. Suppose the contraction map $f_1$ corresponds with this area of $S_C$. Then for any $x$, $y \in S_C$, we compute $|f_1(x) - f_1(y)| = \left(\frac{1}{3}\right)|x - y|$, which means the distance between any $x$ and $y$ in $S_C$ has been reduced by a third. We mimic the rest of the $f_i$'s in the similar way to find each constant to be $\left(\frac{1}{3}\right)$. Hence we use the equation $\sum_{i=1}^{8} C_i^s = 8\left(\frac{1}{3}\right)^s = 1$ to obtain $s = \log 8 / \log 3$.*

If our set $A = \{.0x00xx0000xxxx\ldots : x = 0, 1\}$ is self-similar, then everything above must agree with the properties of our set. But since it is difficult to capture the geometrical behavior of $A$, we move-on to find other means to determine self-similarity.

From [37], if $A$ is self-similar then $d_{bc}(A)$ exists. We know $d_{bc}(A)$ exists whenever $\underline{d}_{bc}(A) = \overline{d}_{bc}(A)$. But we computed two distinct lower and upper box-counting dimensions where the lower box-counting dimension is the correct dimension for our set.

Hence although our set appears to be self-similar, it is not. This is an irregular set whose lower box-counting dimension corresponds with the lower Hausdorff dimension.

# Chapter 7

# Suggestions for Further Investigation

## 7.1 Nonstandard Approach to Hausdorff Measure

In this thesis, our main focus was to introduce a nonstandard approach to the Hausdorff measure so that one may investigate microscopic behavior of measurable functions. Nonstandard methods optimize the coverings that are used so that we obtain a measure of a set as accurate as much as possible.

And now, we ask ourselves some questions regarding the results that which we have obtained which may lead to further studies.

First, under Theorem 5.2.2(a), we said the inequality works only in one direction. So a question that we may once again ask is: are there some sets so that we have an equality in (a) for some appropriate $\delta \approx 0$?

Secondly, we speculate on the uses of Hausdorff dimension in studying connections with box-counting dimension. Because of the particular position of the grids used for Conjecture 5.3.1, it might be necessary that our set must be self-similar when using equation (5.5). If our set is not self-similar, it means that the set must be behaving in some random and unpredictable order. So we must analyze each of the bigger $k$-squares in (5.2). And hence the two-grid approach is avoidable whenever our set is irregular since the classical box-counting method proves to be much easier and faster than equation (5.5).

Thirdly, does equation (5.6) work only for self-similar sets? If there are other sets so that our simplified nonstandard box-counting dimension equals the standard Hausdorff dimension, then what are they? Are we limited to self-similar sets as we have seen under section 5.4.1, or can we extend this to sets with irregular behavior?

We haven't discussed too much of density theory in this thesis but we may also want to approach the Hausdorff density using nonstandard methods because it may give us a better estimate of the bounds for the densities of a set. And since density theory deals with local measure, we might want to look into sets with local measure zero because they might actually be sets with, although small, but positive local measure after using some nonstandard techniques.

We came across a set whose area was an integer with its Hausdorff dimension also being an integer. This was the quadratic Koch curve. Another analysis which we may want to consider is to classify fractals in the order of their Hausdorff dimension using nonstandard analysis. On a separate note, we might ask: is the Hausdorff dimension of some set integral if its length, area, volume, etc is integral?

The equation (5.6) is very similar to the standard box-counting measure. In fact when a set is self-similar, our $d$ equals the Hausdorff dimension. We also believe the residue we obtain might tell us something more about our set, such as, it might give us a slightly better calculation of the dimension. This claim must be further examined.

## 7.2 Returning to Sets Generated by Sequences

There are other sequences that are yet to be explored. Since linear, quadratic, cubic sequences are all of the same class, we're interested in sequences of higher order. So as $k$ vary over the naturals, we have various exponential sequences between $n^k$ and $2^n$. Let $[x]$ denote the integer part of $x$. Then the inequality

$$n < n^2 < \ldots < n^k < \ldots < 2^{[(\log_2 n)(\log_2 \log_2 n)]} < \ldots$$
$$< 2^{\left[(\log_2 n)^2\right]} < c^n < 2^n$$

holds true for $1 < c < 2$.

Now, recall the notation $2_n = 2^{2^{2^{\cdots 2}}}$ where 2 is raised to the power of 2 $(n-1)$-times. Some sequences between $2^n$ and the power-tower of 2 are

$$2^n < 2_{[\log_2{}^{**}n]} < 2_{[\log_2{}^{*}n]} < \ldots$$
$$< 2_{[\log_2 \log_2 n]} < 2_{[\log_2 n]} < c_n < 2_n$$

where $1 < c < 2$. We write $\log_2{}^*n$ to denote the number of times we take $[\log_2(n)]$ to get to 1. Also, $\log_2{}^{**}n$ denotes the number of times we take $\log_2{}^*n$ to equal 1. We may want to explore these intermediate sequences to determine their relationship to the Hausdorff dimension. That is, we might want to look into the changes of the Hausdorff dimension as the sequences change.

Another set that we may dare to investigate has the following conditions:

suppose we have $N \geq 2$ sequences

$$a^1 = \{a_n{}^1\}_{n\in\mathbb{N}}$$
$$a^2 = \{a_n{}^2\}_{n\in\mathbb{N}}$$
$$a^3 = \{a_n{}^3\}_{n\in\mathbb{N}}$$
$$\vdots \quad \vdots$$
$$a^N = \{a_n{}^N\}_{n\in\mathbb{N}}.$$

Sequentially alternate the terms to generate a set $A$ which looks like

$$0.\overbrace{x^1\ldots x^1}^{a_0^1}\overbrace{x^2\ldots x^2}^{a_0^2}\overbrace{x^3\ldots x^3}^{a_0^3}\ldots\overbrace{x^N\ldots x^N}^{a_0^N}\overbrace{x^1\ldots x^1}^{a_1^1}\overbrace{x^2\ldots x^2}^{a_1^2}\ldots\overbrace{x^N\ldots x^N}^{a_1^N}\ldots$$

where each $x^i = 0, 1, \ldots, b_i$ for $0 \leq b_i \in \mathbb{N}$. We haven't thought too much on how we would go about on finding the Hausdorff dimension of $A$, or its self-similarity, or its compactness, but this may be an interesting set to analyze.

Lastly, if $a = \{a_i\}_{i\in\mathbb{N}}$ and $b = \{b_i\}_{i\in\mathbb{N}}$, we may want to examine the results when $a_i \neq b_i$ for almost all $i$. How would the Hausdorff dimension be affected, and would our $\delta$-covering continue to be optimized when we choose it to follow after a chain of zeros, and not the $x$'s? If not, when would we compute the correct Hausdorff dimension, and is there a consistent pattern of choosing the correct $\delta$-coverings as our nondecreasing sequences $a$ and $b$ change at a different rate?

We investigated nonstandard approach to measure theory, and we generated some sets whose Hausdorff dimension is less than 1. Incorporating nonstandard theory to any field of mathematics results in better residue which may lead us to more accurate dimension, measure, etc., or some property that which cannot be observed in the simple standard space. This research does not end here, but is forwarded to all researchers, scientists, students and individuals interested in pure mathematics which will give them a greater appreciation of the world in which we live.